# HJB and Fokker-Planck equations for river environmental management based on stochastic impulse control with discrete and random observation

**Authors** Hidekazu Yoshioka[1, 2], Motoh Tsujimura[3], Kunihiko Hamagami[4], Yuta Yaegashi[5], and Yumi Yoshioka[6]

[1] Assistant Professor, Graduate School of Natural Science and Technology, Shimane University, Nishikawatsu-cho 1060, Matsue, 690-8504, Japan

[2] Center Member, Fisheries Ecosystem Project Center, Shimane University, Nishikawatsu-cho 1060, Matsue, 690-8504, Japan

[3] Professor, Graduate School of Commerce, Doshisha University, Karasuma-Higashi-iru, Imadegawa-dori, Kamigyo-ku, Kyoto, 602-8580, Japan

[4] Associate Professor, Faculty of Agriculture, Iwate University, Morioka, 3-18-8 Ueda, 020-8550, Japan

[5] Independent Researcher, Dr. of Agr., 10-12-403, Maeda-cho, Niihama, 792-0007, Japan

[6] Assistant Professor, Graduate School of Natural Science and Technology, Shimane University, Nishikawatsu-cho 1060, Matsue, 690-8504, Japan

\* Corresponding author
E-mail: yoshih@life.shimane-u.ac.jp

**Abstract**

We formulate a new two-variable river environmental restoration problem based on jump stochastic differential equations (SDEs) governing the sediment storage and nuisance benthic algae population dynamics in a dam-downstream river. Controlling the dynamics is carried out through impulsive sediment replenishment with discrete and random observation/intervention to avoid sediment depletion and thick algae growth. We consider a cost-efficient management problem of the SDEs to achieve the objectives whose resolution reduces to solving a Hamilton-Jacobi-Bellman (HJB) equation. We also consider a Fokker-Planck (FP) equation governing the probability density function of the controlled dynamics. The HJB equation has a discontinuous solution, while the FP equation has a Dirac's delta along boundaries. We show that the value function, the optimized objective function, is governed by the HJB equation in the simplified case and further that a threshold-type control is optimal. We demonstrate that simple numerical schemes can handle these equations. Finally, we numerically analyze the optimal controls and the resulting probability density functions.

## 1. Introduction

The purposes of this paper are to formulate and analyze a cost-efficient environmental restoration problem from a viewpoint of stochastic control with partial observation and partial integro-differential equations. We firstly present the problem background, related issues, and then our contributions.

### 1.1 Problem background

Rivers serve as an essential element on the earth providing habitats for aquatic species, driving hydrological cycles, and providing water resources for human activities [1]. Many rivers are facing with anthropogenic impacts, such as dam construction, industrial pollution, and urbanization, which affect their water quality and quantity, and ecosystems [2-4].

A common issue of rivers having transverse hydraulic structures like dams is the sediment trapping [5-6]; sediment particles flowing along a river are trapped by a dam. The lack of the sediment supply critically affects the downstream river environment. Such examples include the reduction of physical disturbance acting on the riverbed, potentially triggering thick growth of nuisance filamentous benthic algae [7-8] due to the absence of the algae detachment by the sediment transport [9], deterioration of the spawning habitats of fishes [10-11], and artificial river morphological changes [12-13]. The lack of sediment transport has been mitigated by flushing out the stored sediment through operating a dam [14-16] and/or supplying earth and soils from outside the river [17-18]. The timing and amount of replenishment of the sediment potentially depend on the given environmental conditions and should be cost-efficient as well as effective enough for avoiding critical impacts. Therefore, there is a substantial need for establishment of a cost-effective framework of sediment storage management in rivers.

### 1.2 Mathematical background

River environment and ecosystems are driven by stochastic flow discharges and are reasonably considered as stochastic processes [19-22]. Optimal control of system dynamics as stochastic processes can be analyzed in the framework of stochastic control based on the dynamic programming principle [23]. The dynamic programming approach has been widely utilized in modeling and computation of stochastic control problems especially when the dynamics to be optimized follow stochastic differential equations (SDEs): an efficient mathematical tool for describing stochastic phenomena having continuous and/or jump noises [24-27]. An advantage of this approach is that solving a stochastic control problem can reduce to computing a solution governed by the corresponding Hamilton-Jacobi-Bellman (equation) of the partial integro-differential form [23]. A disadvantage of the approach is that often solving the equation must be carried out numerically because of its nonlinear and nonlocal nature. This issue has been overcome by utilizing novel numerical methods in applied problems [28-30]. Theoretical convergence results of numerical methods for HJB and related equations have been compiled in the literature [31-33].

A key property specific to the sediment replenishment problems from the viewpoint of stochastic control is that each replenishment event, transportation of earth and soils from outside the river to the inside it, is considered to be impulsive because the time required for the replenishment (hourly scale) is much shorter than the time scale of the problem (weekly to monthly scales). Therefore, the problem can be formulated as a stochastic impulse control problem [34-35]. There exist many stochastic impulse control models ranging from the exactly solvable [36], semi-analytical [37], and more complex models [38]. Most of them are based on the complete information assumption that the decision-maker who controls the system can continuously receive the complete information of its dynamics.

Another specific point of the sediment replenishment problems, which is common to other environmental and ecological management problems, is that the decision-maker is not always accessible to complete information of the target dynamics, but can collect the information only discretely [39]. Such a problem is a stochastic control problem with partial (discrete and sometimes random) observations being different from the above-mentioned models, and has been used for analyzing epidemic modeling [40], statistical inference [41], insurance management [42-43], and resource management [44].

Recently, based on a stochastic impulse control formalism with Poisson observation times, a single-variable sediment replenishment problem solely focusing on the sediment storage dynamics was formulated [35], and its extension with a replenishment delay has been considered [45]. A severe limitation of these models is that they only consider the sediment storage management, but not the other dynamics like the nuisance benthic algae dynamics [46]. In addition, they assume that the river flow regimes follow Markov chains and that the rate of flushing out of sediment toward downstream is bounded and piecewise constant in each regime. However, this assumption is less consistent with the experimental results that the sediment flushing as well as algae detachment occur impulsively (within few hours) if the river discharge is sufficiently high [47]. This means that the sediment flushing out and the algae detachment occur at the same time as a jump event. In this view, it is physically more reasonable to consider coupled sediment storage and algae population dynamics using jump processes like Lévy processes [48]. To the best of the authors' knowledge, such an attempt has not been reported so far.

**1.3 Objectives and contributions**

The objectives of this paper are to formulate, analyze, and compute a stochastic impulse control problem of a coupled sediment-algae dynamics in a dam-downstream river based on discrete observations. To achieve the objectives, we consider both the 2-D coupled sediment-algae dynamics that will be approached numerically and simplified 1-D sediment dynamics that will be analytically tractable, both of which play complementary roles with each other. As we will explain later, both cases can be formulated in a unified manner. As explained above, we formulate the system dynamics utilizing a jump process representing flushing out of the sediment as well as the associated algae detachment. The dynamics of sediment flushing out and algae detachment are modelled based on our experimental finding and the semi-empirical formula of sediment transport [49]. The dynamics are nonsmooth because of assuming a compact range of the state variables from both physical and biological viewpoints (non-negativity and boundedness of the sediment storage and algae population). The driving jump processes are accordingly set in a state-dependent manner so that the model becomes consistent with this assumption. Notice that the existing jump process models of similar types, especially those in insurance that have been extensively studied, do not employ such an assumption because of handling essentially unbounded processes that are terminated once approaching some ruin boundaries [42-43, 50-51]. Consequently, our model has a unique nonsmooth property.

The decision-maker, which is an environmental manager, seeks for a cost-efficient sediment replenishment policy to avoid sediment depletion as well as thick algae growth. In addition, we focus on a partial observation case where the decision-maker can observe/intervene the system dynamics only discretely and randomly, utilizing a Poisson observation formalism [35]. We formulate an infinite-horizon objective function containing a penalization of the sediment depletion, a penalization of thick algae growth, and a replenishment cost. Finding an optimal control, a replenishment policy minimizing the objective function, reduces to solving an HJB equation that is nonlinear, nonlocal, and nonsmooth. The equation itself seems to be difficult to analyze mathematically but a simplified tractable case is considered to show that the solution to the HJB equation, which is discontinuous at a boundary point representing the sediment depletion, is the value function and that the optimal control is of a threshold type.

The controlled dynamics are also analyzed in this paper. We heuristically derive a Fokker-Planck (FP) equation governing the probability density function (PDF) of the controlled dynamics. A physical consideration of the impulsive replenishment and jump flushing out of the sediment implies that the PDF has Dirac's Deltas along the depletion and full-storage boundaries of the state space. This singularity is verified in the simplified case analytically, and is confirmed numerically by a Monte-Carlo simulation. The FP equations of impulse controls are topics that have not been studied well except for the stochastic resetting models [52] and several limited classes of impulse game and control problems [53-55], and the above-mentioned sediment replenishment models based on the gradual sediment transport [35, 45]. Our FP equation is not found in the literature.

We present simple numerical schemes for computing these partial integro-differential equations. A semi-Lagrangian scheme [56] equipped with a simple but careful evaluation method of the nonlocal jump term [57] is applied to the HJB equation. For the FP equation, a cell-centered upwind finite volume (FV) scheme combined with the conservative jump term treatment [58] and a direct discretization of the singular nature along the boundaries is presented. These schemes are examined against the simplified case, demonstrating that the discontinuity and singularity of the equations are correctly captured. The computational results of the 1-D case demonstrate that discretization of the 1-D case, even if it is the simpler case than the original 2-D case, is not straightforward. Especially, we show that a careful discretization is necessary for the nonlocal terms. We also apply these schemes to a realistic problem where the coefficients and parameters are identified from the experimental and available data. Cost-efficiency of sediment replenishment policies and suppression of the algae growth are concurrently considered using the HJB and FP equations. We compute the FP equation using the finite volume scheme and a Monte-Carlo method separately, which are independent and fundamentally different numerical methods with each other and show a good agreement between them. Limitations of the mathematical and numerical modeling methodologies are also discussed. Consequently, we contribute to formulation, analysis, computation, and application of a new stochastic control model of an engineering problem.

The rest of this paper is organized as follows. The mathematical model is explained and the HJB

and FP equations are presented in Section 2. Simplified cases are considered and verified in this section. The numerical schemes for the two equations are presented and examined in Section 3. Our model is applied to a realistic problem in Section 4. Summary and future perspectives of our research are presented in Section 5. Appendices contain technical parts of this paper.

## 2. Mathematical model
### 2.1 Stochastic differential equations

The SDEs governing the system dynamics are formulated. See, also the conceptual figure (**Figure 1**). We consider a time evolution of sediment storage and nuisance benthic algae population in a downstream reach of a dam blocking sediment transport from the upstream river. Our focus is this 2-D case where both the sediment and algae dynamics are coupled, but a simplified 1-D case considering only the sediment dynamics will also be analyzed in later sections. The formulation of the problem in this and subsequent sub-sections handles the 1-D and 2-D cases in a unified manner. The 1-D case itself is meaningful as it is exactly-solvable and can be used for verifying the numerical schemes as explained above. 1-D sediment replenishment problems are interesting themselves and are important from a viewpoint of sediment storage management and have been studied by the authors in a series of studies [35, 45, 57]. Furthermore, the 1-D case serves as a simplified model where the algae population dynamics are absent. The 1-D problem in our case is seen as a stochastic variant of these previous ones. The 2-D case reduces to the 1-D case by simply omitting the algae dynamics or by setting the initial condition of the algae population defined in the next paragraph to 0.

We assume that the stored sediment is physically flushed out toward downstream by river flows and that the sediment is eventually depleted unless replenished by the decision-maker: an environmental manager. The rate of flushing out, which is called the transport rate here, can be estimated based on physical considerations [49]. The benthic algae are effectively removed from the riverbed by the frictional force exerted by the sediment particles. The detachment is not significant if the sediment is depleted even if the river flow is sufficiently strong [59]. Throughout this paper, we employ the Itô calculus for consistency with the previous research.

The time is denoted as $t \geq 0$. The sediment storage $X = (X_t)_{t \geq 0}$ is valued in $[0, \bar{X}]$. The algae population $Y = (Y_t)_{t \geq 0}$ is valued in $[0, \bar{Y}]$. We set $\bar{X} = \bar{Y} = 1$ assuming a suitable normalization of the variables. Set $D = [0,1]^2$.

Flushing out of the sediment by river flows is assumed to be represented by a pure-jump Lévy process $L = (L_t)_{t \geq 0}$ [60] having the Lévy measure $v$ with only positive jumps:

$$dX_t = -dL_t \text{ for } t > 0, \quad X_0 = x \in [0,1]. \tag{1}$$

We assume the finite activities $0 < \int_{(0,\infty)} v(dz) = \lambda < +\infty$ to simplify the problem. This includes a variety of compound Poisson processes and their finitely many sums [61]. Infinite activities cases may be approximated by the finite activities cases following approximate schemes by compound Poisson processes (Chapter 8 of Cont and Tankov [60]).

This model is still incomplete since the process $X$ possibly becomes negative unless jumps are truncated in some way. Let $\omega_l$ be the $l$ th jump time of $L$ and the jump size at $\omega_l$ be $z_l > 0$. Consider a state-dependent jump process $\hat{L} = (\hat{L}_t)_{t \geq 0}$ having the jump $\min\{z_l, X_{\omega_l-}\}$. In this way, we get $X_{\omega_l} = X_{\omega_l-} - \min\{z_l, X_{\omega_l-}\} = \max\{X_{\omega_l-} - z_l, 0\} \geq 0$ and the non-negativity of $X$ is preserved. We thus replace $L$ by $\hat{L}$ in (1).

The algae growth is based on a deterministic growth effect and a jump decrease effect by the sediment transport (See, also **Section 4** and **Appendix F** for physical explanations.). The SDE of the population dynamics is

$$dY_t = f(Y_t)dt - Y_{t-}d\tilde{L}_t \text{ for } t > 0, \ Y_0 = y \in [0,1], \quad (2)$$

where $\tilde{L}$ is a jump process having the jump size $1 - g(X_{\omega_l-}, z_{\omega_l})$ with $g(X_{\omega_l-}, z_{\omega_l}) = \exp(-\xi \min\{X_{\omega_l-}, z_{\omega_l}\}) \in (0,1]$ where $\xi > 0$, inducing the state jump $Y_{\omega_l} = Y_{\omega_l-} g(X_{\omega_l-}, z_{\omega_l})$, and $f$ represents a generic growth coefficient including the classical logistic model $f(y) = Gy(1-y)$ with $G > 0$, $f(0) = f(1) = 0$, and $f \in C^1([0,1])$ [62]. The coefficient $g$ representing the detachment effect has a "min" operator corresponding to the truncation. The algae detachment under the sediment transport occurs in a relatively short time-scale and is reasonably considered to be an exponential population decrease. In addition, the detachment is not significant in the absence of sediment transport ($X_t = 0$). Therefore, the jump decrease of the algae population using $g$ is physically reasonable. This is essentially a Marcus-type interpretation of multiplicative jumps [63].

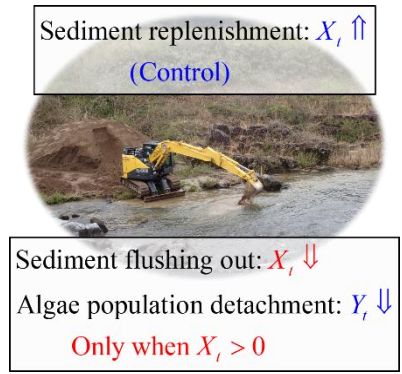

**Figure 1.** A conceptual diagram of our model.

## 2.2 Sediment replenishment

We consider a partial observation/intervention problem where the decision-maker knows the law of

$(X,Y)$ but can observe its sample paths only discretely. See, also **Figure 2**. The simplest model of a discrete observation scheme is the Poisson observation successfully utilized in the conventional models [35, 42-43]. Its advantage is that it has only one and physically interpretable parameter: the observation intensity whose inverse is identified as the mean observation interval. The Poisson nature leads to a model that can be numerically solvable and is exactly-solvable under certain assumptions despite it is nonlinear, nonlocal, and nonsmooth. A Poisson process generating observation times is set as $N = (N_t)_{t \geq 0}$ with the intensity $\Lambda > 0$. The sequence of observation times is $\tau = \{\tau_k\}_{k=0,1,2,...}$, which is increasing with $\tau_0 = 0$.

The decision-maker can collect $(X_{\tau_k}, Y_{\tau_k})$ at each $\tau_k$. A natural filtration generated by the collected information is $\mathcal{F} = (\mathcal{F}_t)_{t \geq 0}$ with $\mathcal{F}_t = \sigma\left\{ (\tau_j, X_{\tau_j}, Y_{\tau_j})_{0 \leq j \leq k}, k = \sup\{j : \tau_j \leq t\} \right\}$. We assume that the decision-maker can supply earth and soils from outside the river at each observation [35]:

$$X_{\tau_k+} = X_{\tau_k} + \eta_k \text{ with } \eta_k = \begin{cases} 0 & (\text{Do nothing}) \\ 1 - X_{\tau_k} & (\text{Fully replenish}) \end{cases} \text{ for } k \geq 1 \text{ and } \eta_0 = 0, \quad (3)$$

where $\eta_k$ represents the amount of sediment supplied at $\tau_k$. Therefore, the decision-maker does nothing or fully replenish the sediment at each observation. The algae population dynamics are only indirectly controlled by the sediment replenishment. The state variables are still living in $D$ by the replenishment because the variables without interventions are clearly non-negative and are not greater than 1 by (3).

Finally, the system of SDEs subject to the sediment replenishment is formulated as

$$d\begin{pmatrix} X_t \\ Y_t \end{pmatrix} = \begin{pmatrix} -d\hat{L}_t \\ f(Y_t)dt - Y_{t-}d\tilde{L}_t \end{pmatrix} + \begin{pmatrix} \bar{\eta}_t dN_t \\ 0 \end{pmatrix} \text{ for } t > 0, \ (X_0, Y_0) \in D, \quad (4)$$

where the process $\bar{\eta} = (\bar{\eta}_t)_{t \geq 0}$ equals $\eta_k$ at $\tau_k$ and equals 0 otherwise. We understand the product term $\bar{\eta}_t dN_t$ at $\tau_k$ as $\eta_k$. A set of admissible control $\mathcal{C}$ contains a continuous-time processes $\bar{\eta} = (\bar{\eta}_t)_{t \geq 0}$ with $\eta_k$ at $\tau_k$ and 0 otherwise, where $\eta_k$ is measurable with respect to $\mathcal{F}_{\tau_k}$.

*Remark 1*
We have $Y_t = 0$ ($t > 0$) if $Y_0 = 0$.

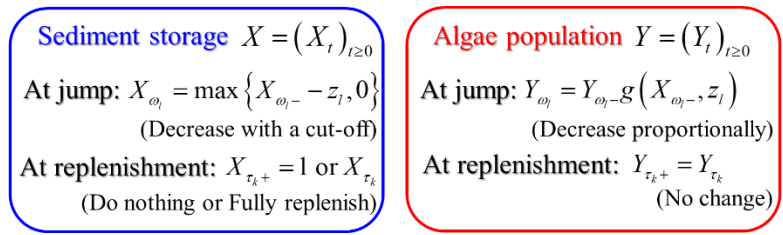

**Figure 2.** A conceptual diagram of the jump and intervention.

## 2.3 Hamilton-Jacobi-Bellman equation

An objective function $\phi: D \times C \to [0, +\infty)$ to be minimized by the replenishment is presented. We consider a long-term management problem and employ an infinite-horizon formulation where both sediment depletion and thick algae growth are concerned. It contains the three terms to be concurrently optimized: a penalization term of the sediment depletion (first term), a penalization term of the thick algae growth (second term), and the replenishment cost (third term):

$$\phi(x, y, \bar{\eta}) = \mathbb{E}^{x,y}\left[\int_0^\infty e^{-\delta s}\chi_{\{X_s=0\}}\mathrm{d}s + \int_0^\infty e^{-\delta s}S(Y_s)\mathrm{d}s + \sum_{k \geq 1} e^{-\delta \tau_k}\left(c\eta_k + d\chi_{\{\eta_k > 0\}}\right)\right], \quad (5)$$

where the conditional expectation with respect to $(X_0, Y_0) = (x, y) \in D$ is denoted as $\mathbb{E}^{x,y}$, $\delta > 0$ is the discount rate, $S: [0,1] \to [0, +\infty)$ is a (uniformly) continuous and increasing function with $S(0) = 0$, $c > 0$ is the coefficient of proportional cost, and $d > 0$ is the fixed cost. Here, $\chi_A$ is a characteristic function of the set $A$: $\chi_A = 1$ if $A$ is true and $\chi_A = 0$ otherwise. A similar cost term has been utilized in a beach nourishment problem by sand filling [64].

The value function $\Phi: D \to [0, +\infty)$ is the minimized $\phi$ with respect to $\bar{\eta} \in C$:

$$\Phi(x, y) = \inf_{\bar{\eta} \in C} \phi(x, y, \bar{\eta}). \quad (6)$$

Clearly, $\Phi$ is non-negative and bounded. In fact, choosing a null-control $\bar{\eta}_0$ with $\eta_i = 0$ ($i \geq 0$) gives

$$0 \leq \Phi(x, y) \leq \phi(x, y, \bar{\eta}_0) \leq \frac{1 + S(1)}{\delta} < +\infty \quad \text{in } D. \quad (7)$$

An optimal policy $\eta^* = \left(\eta_t^*\right)_{t \geq 0} \in C$ minimizing $\phi$ is referred to as the optimal control. The goal of our control problem is to find $\eta^*$ based on the partial observations. By the presence of the fixed cost $d > 0$, any $\eta^*$ should not replenish the sediment if $X_t = 1$ at observations.

### Remark 2

By **Remark 1** and $S(0) = 0$, the control problem is independent from the algae population if $Y_0 = 0$.

We present our HJB equation as a governing partial integro-differential equation of $\Phi$. This equation is derived formally by invoking a dynamic programming principle [35, 45, 65], but can be justified in a simplified case as shown later.

For our model, the HJB equation is formally derived as (**Appendix A**)

$$\delta\Phi - \mathcal{A}\Phi + \Lambda\left(\Phi - \min_{\eta \in \{0, 1-x\}}\left\{\Phi(x+\eta, y) + c\eta + d\chi_{\{\eta > 0\}}\right\}\right) - \chi_{\{x=0\}} - S(y) = 0, \quad (x, y) \in D, \quad (8)$$

where the operator $\mathcal{A}$ is defined for generic sufficiently regular $\Psi = \Psi(x, y)$ as

$$\mathcal{A}\Psi = f(y)\frac{\partial \Psi}{\partial y} - \int_{(0,\infty)}\left\{\Psi(x, y) - \Psi(x - \min\{x, z\}, yg(x, z))\right\}v(\mathrm{d}z). \quad (9)$$

The HJB equation (8) is considered on $D$ since the characteristics of the partial differential terms are not inward and all the terms do not depend on the information outside it. Especially, we can simply omit $f(y)\frac{\partial \Psi}{\partial y}$ along $x = 0, 1$ by $f(0) = f(1) = 0$. The HJB equation has a discontinuous coefficient $\chi_{\{x=0\}}$ and does not have diffusion terms, implying that its solutions may have some singularity such as discontinuity along $x = 0$.

As in the existing models [35, 45, 65], assuming a Markov control $\eta^* = \eta^*\left(X_{\tau_k}, Y_{\tau_k}\right)$ at $\tau_k$, we guess the optimal control

$$\eta^*\left(X_{\tau_k}, Y_{\tau_k}\right) = \arg\min_{\eta \in \{0, 1-X_{\tau_k}\}} \left\{ \Phi\left(X_{\tau_k} + \eta, Y_{\tau_k}\right) + c\eta + d\chi_{\{\eta>0\}} \right\}. \tag{10}$$

In this way, finding an optimal replenishment policy reduces to solving the HJB equation (8). This is verified analytically for a simplified case later.

By invoking the analytical and computational results in the previous mathematical modeling of the sediment replenishment [35, 45], we infer a threshold-type optimal policy

$$\eta^*\left(X_{\tau_k}, Y_{\tau_k}\right) = \begin{cases} 1 - X_{\tau_k} & \left(X_{\tau_k} \leq \overline{x}\left(Y_{\tau_k}\right)\right) \\ 0 & (\text{Otherwise}) \end{cases} \tag{11}$$

with some $\overline{x} : [0,1] \to (-\infty, 1]$ where negative $\overline{x}$ means no replenishment ($\eta^* \equiv 0$). This is an intuitive policy that the sediment should be replenished if the sediment storage is small, and the threshold depends on the algae population. We show that the policy (11) is indeed optimal from both analytical and numerical viewpoints.

### 2.4 Fokker-Planck equation

Analyzing an optimally controlled system dynamics is as important as finding a driving optimal control to deeper understand a stochastic control problem. An FP equation is a governing equation of a PDF of a system of SDEs, and has been extensively analyzed for diffusion processes [66]. However, FP equations for processes driven by state-dependent jumps and/or impulsive controls are still germinating topics [53, 55, 58, 67] possibly because the solutions are highly problem-dependent. A unified formulation for jump-diffusion models for smooth coefficients exist [68]. PDFs may be simply computed using a Monte-Carlo method without considering FP equations; however, one cannot obtain deeper physical insights.

Our problem has the three ingredients, all of which should be reflected in the FP equation. Firstly, both deterministic drift and state-dependent jumps drive the dynamics, which can be handled by considering the existing formalism [58, 67]. Secondly, assuming the control of the threshold type (11), the state would be impulsively reset to a fixed boundary, as in the stochastic resetting models [52]. Finally, related to the first and second ones, the FP equation would not be defined in a classical sense but in a distributional sense along the boundaries $x = 0, 1$ because there will be non-zero resting times of the sediment storage if the storage is depleted or full. This is because both the flushing out and replenishment events follow jump processes with finite activities.

We set the PDF $P = P(t,x,y)$ as a function of $(t, X_t, Y_t) = (t, x, y)$. For each $t \geq 0$, the PDF should satisfy $\int_D P \mathrm{d}x \mathrm{d}y = 1$. By the physical consideration, we assume that $P$ has the form

$$P(t,x,y) = p(t,x,y)\chi_{(0,1)\times[0,1]} + q(t,y)\delta_{\{x=0\}} + r(t,y)\delta_{\{x=1\}}, \quad (12)$$

where the first term represents the PDF in $(0,1)\times[0,1]$, and the second and third terms the distributional terms along $x = 0, 1$, respectively. Here, $\delta_{\{x=x'\}}$ represents the Dirac's Delta concentrated at $x = x'$. The weights $p, q, r$ should be non-negative and integrable. If the PDF is not singular along the boundaries, then simply set $q, r \equiv 0$. The mass conservation property that should be equipped with the PDF is

$$\int_{[0,1]} q(t,y)\mathrm{d}y + \int_{[0,1]} r(t,y)\mathrm{d}y + \int_{(0,1)\times[0,1]} p(t,x,y)\mathrm{d}x\mathrm{d}y = 1, \quad t \geq 0. \quad (13)$$

The FP equation of the controlled dynamics by the threshold-type policy (11) is formulated as

$$\begin{aligned}
&\frac{\partial p}{\partial t} + \frac{\partial}{\partial y}(f(y)p) + \lambda p + \Lambda \chi_{(0,\bar{x}]} p \\
&= \int_{(0,\infty)} \frac{\chi_{\{y \leq g(x,z)\}} \chi_{\{x < 1-z\}}}{g(x,z)} p\left(t, x+z, \frac{y}{g(x,z)}\right) v(\mathrm{d}z), \quad D \cap \{0 < x < 1\}, \\
&\quad + \int_{(0,\infty)} \delta_{\{x=1-z\}} \frac{\chi_{\{y \leq g(x,z)\}}}{g(x,z)} r\left(t, \frac{y}{g(x,z)}\right) v(\mathrm{d}z)
\end{aligned} \quad (14)$$

$$\begin{aligned}
\frac{\partial q}{\partial t} + \frac{\partial}{\partial y}(f(y)q) + \Lambda q &= \int_{(0,1)} \int_{(0,\infty)} \chi_{\{x \leq z\}} p(t,x,y) v(\mathrm{d}z) \mathrm{d}x \\
&\quad + \int_{(0,\infty)} \frac{\chi_{\{y \leq g(x,z)\}} \chi_{\{z \geq 1\}}}{g(1,z)} r\left(t, \frac{y}{g(1,z)}\right) v(\mathrm{d}z)
\end{aligned}, \quad D \cap \{x = 0\}, \quad (15)$$

$$\frac{\partial r}{\partial t} + \frac{\partial}{\partial y}(f(y)r) + \lambda r - m\Lambda = 0, \quad D \cap \{x = 1\}, \quad (16)$$

where

$$m = q(t,y) + \int_{(0,\bar{x}]\times[0,1]} p(t,x,y)\mathrm{d}x\mathrm{d}y. \quad (17)$$

The dependence of $\bar{x}$ on $y$ is suppressed here, and we formally omit the terms multiplied by $\Lambda$ and $m = 0$ if $\bar{x}(y) < 0$ (No replenishment). The equations (14)-(17) possess clear physical meanings as explained below. The temporal differential terms simply represent the rate of the changes of the weights with respect to the time. The first-order partial differential terms with respect to $y$ represent the deterministic drift of the PDF by the algae population growth.

The terms multiplied by $\Lambda$ represent the sediment replenishment. The last terms in the left-hand sides of (14) and (15) represents the replenishment from the current state if $x \in [0, \bar{x}]$. The last term in the left-hand side of (17) corresponds to the replenishment from the state just before, where the coefficient $m$ in (17) means the fact that the replenishment is possible only when $x \in [0, \bar{x}]$. These terms on the replenishment are just an application of the stochastic resetting formalism [52].

The remaining terms represent jump decreases of the sediment storage and algae population. The terms multiplied by $\lambda$ simply represents the jump decrease of the state variables at the mean intensity $\lambda$. Clearly, there is no decrease of the sediment if $x = 0$ as shown in (15). The integral terms represent the jump inflow from the state just before to the current state. The appearances of the characteristic functions are due to the compactness of the domain $D$. For example, $\chi_{\{y \leq g(x,z)\}}$ means that the current state $(x, y)$ can only be realized by jumps such that $y \leq g(x, z)$. Indeed, for the jump size $z$, the state just before the current time should be $\left( x+z, \dfrac{y}{g(x,z)} \right) \in D$. The functional forms of these terms can be derived by considering the physics-based consideration considering both pre- and post-jump states (Appendix 2 of Yoshioka et al. [58]).

The following proposition proves that the FP equation is indeed conservative, meaning that each term is correctly balanced. Although this property seems to be trivial for conventional FP equations, it is not so for our case because of its unique form. Its proof is placed in **Appendix B**.

*Proposition 1*
*If the equality in (13) is satisfied at $t = 0$, then it holds true for $t > 0$.*

### 2.5 Analysis of a simplified case

We presented the HJB and FP equations in the previous sub-sections. We consider a simplified problem justifying them as well as demonstrating their singular nature. It seems to be impossible to obtain exact solutions to the HJB equations and FP equations for the genuinely 2-D case; however, we can exactly solve them if we focus on the no-algae case $Y_0 = y = 0$. In this case, clearly $Y \equiv 0$ and the problem simply becomes a sediment replenishment problem to avoid the depletion. We show that the singular nature of the HJB and FP equations emerge even in this case. The reduced equations are not found in the literature.

In the 1-D case, the objective function becomes

$$\phi(x, \bar{\eta}) = \mathbb{E}^x \left[ \int_0^\infty e^{-\delta s} \chi_{\{X_s = 0\}} \mathrm{d}s + \sum_{k \geq 1} e^{-\delta \tau_k} \left( c\eta_k + d \chi_{\{\eta_k > 0\}} \right) \right], \tag{18}$$

and the corresponding value function as a map from $[0,1]$ to $[0, +\infty)$ is

$$\Phi(x) = \inf_{\bar{\eta} \in \mathcal{C}} \phi(x, \bar{\eta}), \tag{19}$$

where the conditional expectation with respect to $X_0 \in [0,1]$ is denoted as $\mathbb{E}^x$.

Assume that the controlled dynamics eventually approach to a stationary state. The HJB and FP equations under the simplification are

$$\delta \Phi + \int_{(0,\infty)} \left\{ \Phi(x) - \Phi(x - \min\{x, z\}) \right\} v(\mathrm{d}z) + \Lambda \left( \Phi - \min_{\eta \in \{0, 1-x\}} \left\{ \Phi(x+\eta) + c\eta + d\chi_{\{\eta > 0\}} \right\} \right) - \chi_{\{x = 0\}} = 0 \tag{20}$$

and

$$\lambda p + \Lambda \chi_{(0,\bar{x}]} p = \int_{(0,\infty)} \chi_{\{x<1-z\}} p(x+z) v(dz) + \left( \int_{(0,\infty)} \delta_{\{x=1-z\}} v(dz) \right) r, \qquad (21)$$

$$\Lambda q = \int_{(0,1)} \int_{(0,\infty)} \chi_{\{x\le z\}} p(x) v(dz) dx + \left( \int_{(0,\infty)} \chi_{\{z\ge 1\}} v(dz) \right) r, \quad \lambda r - m\Lambda = 0, \qquad (22)$$

where the dependence of the variables on $y$ is omitted.

We found that the equations (20) and (21)-(22) are still not solvable analytically, and that we can get the exact solutions to these equations if we set $v(dz) = \lambda \chi_{(0,1)} dz$. This assumption seems to simplify the problem too much, but the simplified model shares the expected singularity of the original 2-D model as shown in the propositions below. Especially, **Proposition 2** verifies the optimality of the control (11). Proofs of **Propositions 2-3** are in **Appendices C-D**. For the assumption of **Proposition 2**, see **Remark 3**. Hereafter, for any functions $F$ having the range $[0,1]$, $F(+0)$ (resp., $F(1-)$) is the right-limit of $F$ at $x = 0$ (resp., the left-limit of $F$ at $x = 1$).

*Proposition 2*

Assume $v(dz) = \lambda \chi_{(0,1)} dz$. Then, the bounded function $\bar{\Phi} : [0,1] \to \mathbb{R}$ with $\bar{\Phi} \in C(0,1]$ of the form

$$\bar{\Phi}(x) = \begin{cases} \Phi_0 & (x = 0) \\ \left( \Phi_{+0} - \Phi_0 - \dfrac{c}{\lambda} \Lambda \right) e^{\gamma x} + \Phi_0 + \dfrac{c}{\lambda} \Lambda & (0 < x \le \bar{x}) \\ -(\Phi_0 - \Phi_1) e^{\beta(x-1)} + \Phi_0 & (\bar{x} < x \le 1) \end{cases} \qquad (23)$$

satisfies the HJB equation (20) if there is a unique quintet $(\bar{\Phi}(0), \bar{\Phi}(+0), \bar{\Phi}(1), \bar{x}) = (\Phi_0, \Phi_{+0}, \Phi_1, \bar{x})$ with $\bar{x} \in (0,1)$ satisfying

$$\begin{aligned}
(\delta + \lambda + \Lambda) \Phi_{+0} - \lambda \Phi_0 &= \Lambda \{\Phi_1 + c + d\} \\
(\delta + \Lambda) \Phi_0 &= \Lambda \{\Phi_1 + c + d\} + 1 \\
\left( \Phi_{+0} - \Phi_0 - \dfrac{c}{\lambda} \Lambda \right) e^{\gamma \bar{x}} + \Phi_0 - \Phi_1 + \dfrac{c}{\lambda} \Lambda &= c(1 - \bar{x}) + d \quad \text{with } \beta = \lambda/(\delta + \lambda), \ \gamma = \lambda/(\delta + \lambda + \Lambda). \ (24) \\
(\Phi_0 - \Phi_1)\left( 1 - e^{\beta(\bar{x}-1)} \right) &= c(1 - \bar{x}) + d
\end{aligned}$$

*Furthermore, we have $\bar{\Phi} = \Phi$ and that the control (11) is optimal.*

*Proposition 3*

Assume $v(dz) = \lambda \chi_{(0,1)} dz$ and $0 \le \bar{x} < 1$. If $0 < \bar{x} < 1$, then the FP equation (21)-(22) has the solution:

$$r = \left( \lambda \Lambda^{-1} + e^{1-\bar{x}} \right)^{-1} > 0, \quad q = r \left( \lambda \Lambda^{-1} - e^{1-\bar{x}} \left( e^{\alpha \bar{x}} - 1 \right) \right) > 0, \qquad (25)$$

$$p(x) = \begin{cases} \alpha r e^{1-x+\alpha \bar{x} - \alpha x} & (0 < x \le \bar{x}) \\ r e^{1-x} & (\bar{x} < x < 1) \end{cases} \quad \text{with } \alpha = \lambda/(\lambda + \Lambda). \qquad (26)$$

*If $\bar{x} = 0$, then we can simply take the limit $\bar{x} \to +0$ in (25) and omit the first line of (26). For $0 \leq \bar{x} < 1$, this solution is discontinuous at $x = \bar{x}$ and has Dirac's Deltas at $x = 0,1$.*

*Remark 3*

The system (24) in **Proposition 2** reduces to a single equation of $\bar{x}$:

$$(c+d-c\bar{x})l(x) = \left(\frac{1}{\delta+\lambda+\Lambda} + \frac{\Lambda}{\lambda}c\right)e^{\gamma\bar{x}} - \frac{\Lambda}{\lambda}c \quad \text{with} \quad l(x) = \frac{e^{-\beta(1-\bar{x})}}{1-e^{-\beta(1-\bar{x})}}. \tag{27}$$

The left- and right-hand sides of (27) are denoted as $f_L(\bar{x})$ and $f_R(\bar{x})$ with $\bar{x} \in [0,1)$. We get $f_L(1-) = +\infty$, $0 < f_R(1-) < +\infty$, $f_L(0) = (c+d)\frac{e^{-\beta}}{1-e^{-\beta}} > 0$, and $f_R(0) = \frac{1}{\delta+\lambda+\Lambda} > 0$. Therefore, for small $c+d > 0$, we have $f_L(0) < f_R(0)$ and $f_L(1-) > f_R(1-)$, and get existence of $\bar{x} \in (0,1)$ solving (27) by classical intermediate value theorem. By $l' = \beta l(1+l)$, we see

$$f'_L = l(-c + \beta(1+l)(c+d-cx)) \quad \text{and} \quad f''_L = \beta l(1+l)(-2c + \beta(c+d-cx)(1+2l)). \tag{28}$$

For small $c > 0$, by (28) and $l', l > 0$, we get $f'_L, f''_L > 0$ in $[0,1)$. This is because we can always choose small $c > 0$ and $c/d$ such that $-c + \beta d > 0$ and $-2c + \beta d > 0$. In summary, if $c+d$ and $c/d$ are sufficiently small, then (27) has a unique solution $\bar{x} \in (0,1)$. The system (24) is then uniquely solvable because the other unknowns are obtained using this $\bar{x}$.

**Proposition 2** shows that the optimal policy is of the threshold type under the unique solvability assumption of the system of (24). This is a non-trivial issue, but can be at least numerically analyzed using some existing methodologies like a Picard method, Newton method, bisection method, or a trial and error method. We do not analyze this issue deeper, instead numerically demonstrate that the system can be handled accurately. The verification result would apply to the 2-D case if we assume a uniformly bounded solution to the HJB equation having a discontinuity only along $x = 0$. Its proof is essentially the same with that in **Appendix B** except for that each term becomes a bit complicated.

**Proposition 3** shows that the FP equation on the other hand is completely exactly solvable given a threshold $\bar{x}$. The solution is piecewise continuous and is discontinuous at the threshold level, and admits Dirac's Deltas at the boundaries as. This implies that solving the FP equation, even in this simplified case, requires accurately approximating the singular nature. In this paper, we therefore propose a numerical scheme that discretize the PDF in a cell-centered manner inside the domain, while the boundary singularity is handled by directly solving the equations (15)-(16).

**Figures 3** and **4** plot the exact solutions to the HJB and FP equations under the simplifications. The parameter values are set as $\delta = 0.1$, $c = 0.35$, $d = 0.30$, $\Lambda = 0.25$, $\lambda = 0.20$. These values have been chosen so that the discontinuities and singularities of the solutions to the HJB and FP equations, which are their characteristics as well as potential difficulties when numerically computing them, are visible. We get $\bar{x} = 0.7986$ using a bisection method with the error smaller than $10^{-13}$. We get $\Phi_0 = 4.253$, $\Phi_{+0} = 2.435$, and $\Phi_1 = 1.304$. For the FP equation, the corresponding weights of the Dirac's Deltas are $q = 0.138$ and $r = 0.494$. A Monte-Carlo result using a forward-Euler method (time increment 0.02) with $6 \times 10^6$ sample paths is also plotted in **Figure 4**, verifying the derived exact PDF inside the domain. The computed $q$ and $r$ are 0.137 and 0.495, respectively, again verifying the exact one.

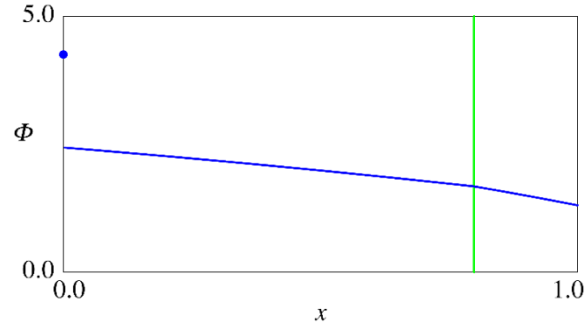

**Figure 3.** The exact value function of the HJB equation in the reduced case (Blue) and the threshold (Green).

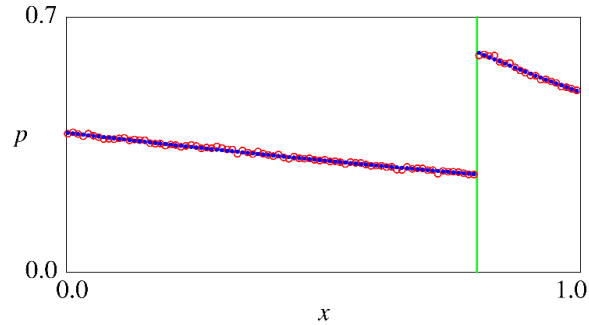

**Figure 4.** The PDFs of the reduced FP equation inside the domain (Blue: exact solution, Red: Monte-Carlo simulation) and the threshold (Green).

## 3. Numerical schemes

Numerical schemes for discretizing the HJB and FP equations are presented. The underlying techniques are not new, but the equations to be discretized are new. See **Figure 5** for the computational domain.

### 3.1 Numerical scheme for the HJB equation

The numerical scheme for computing the HJB equation is based on the semi-Lagrangian scheme with piecewise linear interpolations equipped with the Weighted Essentially Non-Oscillatory reconstruction (WENO) and the conventional value iteration [56]. The nonlocal term corresponding to jumps is discretized using piecewise constant treatment [57] that can handle the boundary discontinuity for simpler problems. Semi-Lagrangian schemes have been widely applied to degenerate elliptic/parabolic equations arising in stochastic control and related problems [69-70].

Since the discretization of the semi-Lagrangian part is essentially the same with Carlini et al. [56] and that for the "min" part with Yoshioka and Yaegashi [45], the explanation here focuses on discretization of the nonlocal term:

$$\mathcal{N}\Phi = -\int_{(0,\infty)} \{\Phi(x,y) - \Phi(x - \min\{x,z\}, yg(x,z))\} v(\mathrm{d}z). \tag{29}$$

Set the vertices $P_{i,j}:(x_i, y_j) = (ih, jh)$ ($1 \leq i, j \leq n$) with $n \in \mathbb{N}$ and $h = n^{-1}$. Set $x_i = ih$ and $y_j = jh$ ($0 \leq i, j \leq n$). Set $\bar{y}_j = (j - 0.5)h$ ($1 \leq i, j \leq n$). The space of jump size is approximated by a bounded set $Z = (\underline{z}, \bar{z})$ with some $0 \leq \underline{z} < \bar{z} < +\infty$. If the jump size is unbounded, $\bar{z}$ is chosen sufficiently large. The measure $v$ should be approximated as well if necessary. The space $Z$ is discretized into $L \in \mathbb{N}$ intervals $Z_l = (\underline{z} + (\bar{z} - \underline{z})(l-1)L^{-1}, \underline{z} + (\bar{z} - \underline{z})lL^{-1})$ ($1 \leq l \leq L$). Set $z_l = \underline{z} + (\bar{z} - \underline{z})(l - 0.5)L^{-1}$ ($1 \leq l \leq L$) and $v_l = \int_{Z_l} v(\mathrm{d}z)$ with $\sum_{l=1}^{L} v_l = \lambda$. The solution $\Phi$ is discretized at each $P_{i,j}$, and the approximation at this point is denoted as $\Phi_{i,j}$. Invoking Yoshioka and Tsujimura [57], we discretize (29) as

$$(\mathcal{N}\Phi)_{i,j} = -\lambda \Phi_{i,j} + \sum_{l=1}^{L} v_l \Phi_{a_{i,l}, b_{i,j,l}} \tag{30}$$

with $a_{i,l} = \lfloor n(x_i - z_l) \rfloor$ and $\omega_{j,l} = \lfloor n\bar{y}_j g(x_i, z_l) \rfloor$, where $\lfloor \cdot \rfloor$ represents the integer that is equal to or not larger than the argument.

The other nonlocal term representing the replenishment is discretized as

$$(\mathcal{R}\Phi)_{i,j} = \Lambda \left( \Phi_{i,j} - \min_{\hat{\eta} \in \{0, n-i\}} \left\{ \Phi_{i+\hat{\eta}, j} + c\hat{\eta}h + d\chi_{\{\hat{\eta} > 0\}} \right\} \right). \tag{31}$$

By (31), the optimal control at $P_{i,j}$, denoted as $\eta^*_{i,j}$, is

$$\eta^*_{i,j} = h \times \arg\min_{\hat{\eta} \in \{0, n-i\}} \left\{ \Phi_{i+\hat{\eta}, j} + c\hat{\eta}h + d\chi_{\{\hat{\eta} > 0\}} \right\}. \tag{32}$$

For each $j \geq 0$ if there is some $i' \geq 0$ with $\eta^*_{i',j} = 1$ ($0 \leq i \leq i'$) and $\eta^*_{i',j} = 0$ ($i' + 1 \leq i \leq n$), then we detect the free boundary $\bar{x}_j = 0.5(x_j + x_{j+1})$. Otherwise, set $\bar{x}_j = -1$. The sequence $\{\bar{x}_j\}_{0 \leq j \leq n}$ is used in the computation of the FP equation.

Using the pseudo-time increment $\rho > 0$ for the semi-Lagrangian discretization, the HJB equation at each $P_{i,j}$ ($0 \leq i, j \leq n$) is discretized as

$$\Phi_{i,j} = (\mathcal{V}\Phi)_{i,j} \equiv e^{-\delta\rho}(\Pi\Phi)_{i,j} + \frac{1 - e^{-\delta\rho}}{\delta} \left\{ (\mathcal{N}\Phi)_{i,j} - (\mathcal{R}\Phi)_{i,j} + \chi_{\{i=0\}} + S(y_j) \right\}, \tag{33}$$

where $(\Pi\Phi)_{i,j}$ is the approximation of $(\Phi(x, y + f(y)\rho))_{i,j}$ evaluated with the semi-Lagrangian scheme using the WENO reconstruction in the $y$-direction [56]. The value iteration is then applied to (33) as $\Phi_{i,j}^{(K)} = (\mathcal{V}\Phi)_{i,j}^{(K)}$ ($K = 0,1,2,...$) starting from some initial guess, which is $\Phi_{i,j}^{(0)} \equiv 0$ in this paper. The numerical parameter $\rho$ is chosen as $h^{1.5}$ to achieve at least the first-order accuracy [56].

### 3.2 Numerical scheme for the FP equation

Our FP equation is singular along $x = 0,1$. This characteristic is directly considered in the numerical scheme based on (14)-(17). We solve (14) using a cell-centered FV scheme that directly handles the equations (14)-(17) inside the domain and along the boundaries separately.

The $(i,j)$ th cell $C_{i,j}$ ($1 \le i, j \le n$) is the open rectangle surrounded by the four vertices $P_{i-1,j-1}$, $P_{i,j-1}$, $P_{i,j}$, and $P_{i-1,j}$. These cells are utilized for discretizing $p$, and that approximated on $C_{i,j}$ is denoted as $p_{i,j}$. The boundaries $x = 0,1$ are also discretized to approximate $q, r$. The interval $I_j^{(L)}$ ($1 \le j \le n$) (resp., $I_j^{(R)}$ ($1 \le j \le n$)) is defined as the open interval bounded by the two vertices $P_{0,j-1}$ and $P_{0,j}$ (resp., $P_{n,j-1}$ and $P_{n,j}$). The discretized $q$ on $I_j^{(L)}$ ($r$ on $I_j^{(R)}$) is denoted as $q_j$ (resp., $r_j$). The time increment for temporal discretization is $\Delta > 0$.

Set the numerical fluxes $F_{p,i,j} = f(\bar{y}_j)\hat{p}_{i,j}$ ($1 \le i \le n$, $1 \le j \le n$), $F_{q,i,j} = f(\bar{y}_j)\hat{q}_j$ and $F_{r,i,j} = f(\bar{y}_j)\hat{r}_j$ ($1 \le j \le n$). They are set to be 0 if $j = 0, n$. The quantities $\hat{p}_{i,j}, \hat{q}_j, \hat{r}_j$ are evaluated by the classical first-order upwind method with the WENO reconstruction [71]. Approximate $\bar{x}_j$ (approximated one in the previous section or an exact one if available) by $\theta_j = \lfloor \bar{x}_j n \rfloor h$. The difference between $\bar{x}_j$ and $\theta_j$ is not larger than $h$.

The spatial semi-discretization of the FV scheme on $C_{i,j}$ is

$$\frac{\mathrm{d}p_{i,j}}{\mathrm{d}t} + \Lambda \chi_{\{i \le \theta_j\}} p_{i,j} + \frac{1}{h}(F_{p,i,j} - F_{p,i,j-1}) + \lambda p_{i,j} - J_{i,j} = 0 \qquad (34)$$

and those on $I_j^{(L)}$ and $I_j^{(R)}$ are

$$\frac{\mathrm{d}q_j}{\mathrm{d}t} + \Lambda q_j + \frac{1}{h}(F_{q,j} - F_{q,j-1}) - J_j^{(L)} = 0, \qquad (35)$$

$$\frac{\mathrm{d}r_j}{\mathrm{d}t} + \frac{1}{h}(F_{r,j} - F_{r,j-1}) + \lambda r_j - \Lambda m_j = 0 \quad \text{with} \quad m_j = q_j + \sum_{i=1}^{\theta_j} p_{i,j} h, \qquad (36)$$

respectively. The terms $J_{i,j}$ and $J_j^{(L)}$ are the discretized integrals representing the jump inflows. Their discretization is a straightforward extension of Yoshioka et al. [58], which is similar to that of the conservative semi-Lagrangian scheme using piecewise constant basis [72]. Set $\alpha_{i,l} = \lfloor n(\bar{x}_i - z_l) \rfloor$, $\beta_{i,j,l} = \lfloor n\bar{y}_j g(\bar{x}_i, z_l) \rfloor$, $\gamma_l = \lfloor n(1 - z_l) \rfloor$, and $\omega_{j,l} = \lfloor n\bar{y}_j g(1, z_l) \rfloor$, which are integers to identify the cells after each jump from the cell $C_{i,j}$. The two cases appear: $\alpha_{i,l} \ge 0$ and $\alpha_{i,l} < 0$. The former is a jump to inside $D$, while the latter to the boundary $x=0$. Similarly, we have the cases $\gamma_l \geq 0$ and $\gamma_l < 0$. We discretize each jump term as follows:

$$J_{i,j} = \sum_{l=1}^{L}\sum_{i',j'=1}^{n} \chi_{\{\alpha_{i',l}=i-1,\beta_{i',j',l}=j-1\}} v_l p_{i',j'} + \sum_{l=1}^{L}\sum_{j'=1}^{n} \chi_{\{\gamma_l=i-1,\omega_{j',l}=j-1\}} v_l r_{j'} h^{-1}, \qquad (37)$$

$$J_j^{(L)} = \sum_{l=1}^{L}\sum_{i',j'=1}^{n} \chi_{\{\alpha_{i',l}<0,\beta_{i',j',l}=j-1\}} p_{i',j'} h + \sum_{l=1}^{L}\sum_{j'=1}^{n} \chi_{\{\gamma_l<0,\omega_{j',l}=j-1\}} r_{j'}. \qquad (38)$$

These terms have been constructed considering points of the pre- and post-jumps [58]. Substituting (37) and (38) into (34) and (35) completes the spatial semi-discretization.

The resulting semi-discretized system is fully discretized using the classical forward-Euler method with a sufficiently small time increment for stable computation. The initial condition should satisfy

$$M = \sum_{j=1}^{n} q_j h + \sum_{i,j=1}^{n} p_{i,j} h^2 + \sum_{j=1}^{n} r_j h = 1 \quad \text{at} \quad t=0. \qquad (39)$$

We can show that, as in the continuous counterpart, the semi-discretized FP equation is conservative (See, **Appendix E** for its proof). Like **Proposition 1**, this property seems to be non-trivial because our case because of its unique form. Its classical forward-Euler counterpart is conservative as well.

*Proposition 4.*

*If the equality in (39) is satisfied at $t=0$, then it holds true for $t>0$.*

*Remark 4*

As discussed in **Remark 1**, the population dynamics would be restricted along $y=0$ if the population is initially absent. In this case, the FP equation would be singular along this boundary as well. An option to handle this kind of problems is to consider an additional Dirac Delta defined along this boundary. However, it is not carried out in this paper because our focus is the coupled problem.

*Remark 5*

The algae population may eventually converge to the boundary $y=1$, but such an event does not occur at any finite times. The boundaries $x=0,1$, on the other hand, can be attainable within finite times. Therefore, the boundaries $y=1$ and $x=0,1$ are qualitatively different with each other in this sense.

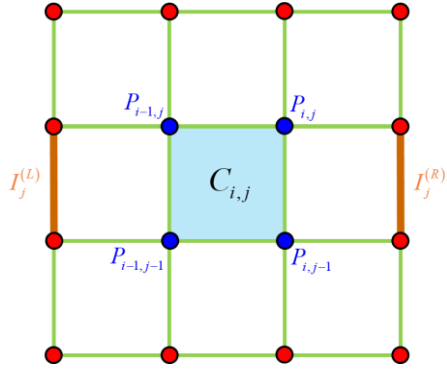

**Figure 5.** A conceptual diagram of the numerical discretization.

## 4. Numerical computation

### 4.1 Convergence analysis

Computational accuracy of the two numerical schemes is analyzed using the exact solutions plotted in **Figures 3** and **4**. We examine $n = 50j$ ( $j = 1, 2, 4, 8, 16$ ) and compute the conventional $l^1$, $l^2$, and $l^\infty$ error norms for each computational resolution. Intuitively, the space of the jump size should be refined as the spatial resolution becomes finer. We set $L = 2n = 100j$, and later examine coarser discretization methods of the jump size, implying that the schemes are not unconditionally convergent. The parameters for the (pseudo-)time discretization are set as $\rho = h = n^{-1}$ and $\Delta = 2n^{-1} = 2h$. The value iteration for the HJB equation is terminated when we get $\left|\Phi^{(K)} - \Phi^{(K-1)}\right| < 10^{-12}$ at all the vertices. The FP equation is temporally integrated in a sufficiently long time because the interest here is the stationary PDF. The temporal integration is terminated if $\left|\phi^{(k)} - \phi^{(k-1)}\right| < 10^{-10}$ for $\phi = p$, $q$, and $r$ at all the vertices.

**Table 1** shows the error norms and convergence rates for the numerical solutions to the HJB and FP equations, respectively. The computational results imply that the numerical scheme for the HJB equation is first-order. The numerical solutions do not have spurious oscillations. The first-order convergence of the HJB equation is consistent with the fact that the scheme is using a piece-wise constant treatment of the jump term. The convergence rate is not so different from conventional schemes for HJB equations having continuous solutions [73-74].

For the FP equation, from **Table 2**, we infer that the convergence rates are first-order in the $l^1$ error norm, 0.5-th order in the $l^2$ error norm, and non-convergent in the $l^\infty$ error norm possibly due to the lower regularity of the PDF than the value function. In fact, for the reference solutions, the value function is discontinuous but bounded, while the PDF is discontinuous and singular. Nevertheless, the convergence of the FV scheme is considered to be satisfactory because the PDF is integrable, and more specifically, it is mass-conservative, but it is not necessarily square-integrable. Note that some super-convergent phenomenon is observed at the resolution $n = 400$, but it is optimistic because the convergence rates are smaller in the other cases. In fact, the estimated CRs using the results with

$n = 200$, $800$, and $1600$ are 1.0 and 0.5 for $l^1$ and $l^2$ error norms, respectively. The reason of this super-convergence may be due to the relatively simple functional forms of the coefficient.

Computational accuracy of the threshold $\bar{x}$ is also important because it determines the optimal threshold policy. **Table 5** presents the computed $\bar{x}$. The exact $\bar{x}$ is 0.7986. The obtained results suggest that the present finite difference scheme for the HJB equation can capture the threshold $\bar{x}$ with an error smaller than the cell size.

The impacts of the discretization of the space of jump size are also discussed. **Tables 3** and **4** show the error norms for the numerical solutions to the HJB and FP equations with the different values of $L$ while the spatial resolution is fixed to $n = 200$. The results demonstrate that $L$ should be chosen not smaller than $n$ and $2n$ for the HJB and FP equations, respectively. This difference between the two equations is considered to be due to using different approximation basis. As shown in **Figures 6** and **7**, using a too small $L$ leads to a slightly stair-like numerical solutions for the HJB equation, while spuriously oscillating solutions to the FP equation. Nevertheless, the threshold is reasonably captured in the HJB equation with the under-resolved jump space discretization (**Figure 6**). The computational results support the resolution $L = 2n$ used in **Tables 1** and **2**.

In summary, the computational results demonstrate that the presented numerical schemes can handle the HJB and FP equations.

**Table 1.** The error norms for the HJB equation with respect to the different values of $n$. Convergence rates (CRs) are also presented in the table. CR for the discretization level $n$ is $\log_2(e_n/e_{n+1})$ with $e_n$ the corresponding error at the level $n$. The CRs in the other tables are calculated in the same way.

| $n$ | 50 | 100 | 200 | 400 | 800 | 1600 |
|---|---|---|---|---|---|---|
| $l^1$ error | 1.383.E-02 | 6.891.E-03 | 3.442.E-03 | 1.720.E-03 | 8.597.E-04 | 4.310.E-04 |
| $l^2$ error | 1.388.E-02 | 6.916.E-03 | 3.454.E-03 | 1.726.E-03 | 8.628.E-04 | 4.326.E-04 |
| $l^\infty$ error | 1.680.E-02 | 8.370.E-03 | 4.180.E-03 | 2.090.E-03 | 1.050.E-03 | 5.300.E-04 |
| $l^1$ CR | 1.01 | 1.00 | 1.00 | 1.00 | 1.00 | |
| $l^2$ CR | 1.01 | 1.00 | 1.00 | 1.00 | 1.00 | |
| $l^\infty$ CR | 1.01 | 1.00 | 1.00 | 0.99 | 0.99 | |

**Table 2.** The error norms for the FP equation with respect to the different values of $n$. Convergence rates (CRs) are also presented in the table.

| $n$ | 50 | 100 | 200 | 400 | 800 | 1600 |
|---|---|---|---|---|---|---|
| $l^1$ error | 5.318.E-03 | 2.656.E-03 | 1.945.E-03 | 8.804.E-05 | 5.060.E-04 | 2.473.E-04 |
| $l^2$ error | 3.032.E-02 | 2.170.E-02 | 2.358.E-02 | 1.652.E-04 | 1.186.E-02 | 8.390.E-03 |
| $l^\infty$ error | 2.182.E-01 | 2.189.E-01 | 3.351.E-01 | 4.015.E-04 | 3.357.E-01 | 3.358.E-01 |
| $l^1$ CR | 1.00 | 0.45 | 4.47 | -2.52 | 1.03 | |
| $l^2$ CR | 0.48 | -0.12 | 7.16 | -6.17 | 0.50 | |
| $l^\infty$ CR | 0.00 | -0.61 | 9.71 | -9.71 | 0.00 | |

**Table 3.** The error norms for the HJB equation with respect to the different values of $L$. The spatial resolution is fixed to $n=200$.

| $L$ | $n/4$ | $n/2$ | $n$ | $2n$ | $4n$ |
|---|---|---|---|---|---|
| $l^1$ | 1.184.E-02 | 8.639.E-03 | 3.441.E-03 | 3.442.E-03 | 3.441.E-03 |
| $l^2$ | 1.259.E-02 | 9.036.E-03 | 3.454.E-03 | 3.454.E-03 | 3.454.E-03 |
| $l^\infty$ | 2.034.E-02 | 1.655.E-02 | 4.180.E-03 | 4.180.E-03 | 4.180.E-03 |

**Table 4.** The error norms for the FP equation with respect to the different values of $L$. The spatial resolution is fixed to $n=200$.

| $L$ | $n/4$ | $n/2$ | $n$ | $2n$ | $4n$ |
|---|---|---|---|---|---|
| $l^1$ | 3.831.E-01 | 7.935.E-03 | 8.176.E-03 | 1.945.E-03 | 1.945.E-03 |
| $l^2$ | 4.946.E-01 | 4.630.E-02 | 2.595.E-02 | 2.358.E-02 | 2.358.E-02 |
| $l^\infty$ | 1.479.E+00 | 3.333.E-01 | 3.367.E-01 | 3.351.E-01 | 3.351.E-01 |

**Table 5.** Computed threshold $\bar{x}$ and the absolute errors between the computed and exact ($\bar{x}=0.7986$) ones.

| $n$ | 50 | 100 | 200 | 400 | 800 | 1600 |
|---|---|---|---|---|---|---|
| Cell size | 0.02 | 0.01 | 0.005 | 0.0025 | 0.00125 | 0.000625 |
| Computed $\bar{x}$ | 0.8100 | 0.7950 | 0.7975 | 0.7988 | 0.7981 | 0.7984 |
| Error | 0.0114 | 0.0036 | 0.0011 | 0.0001 | 0.0005 | 0.0002 |

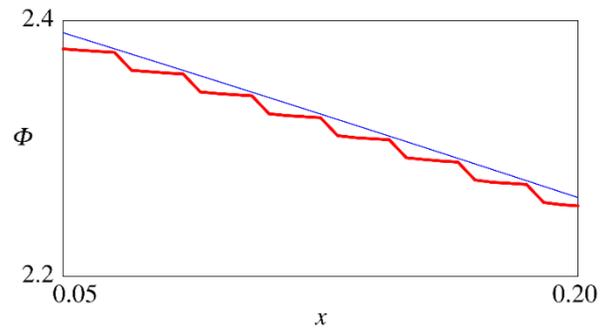

**Figure 6.** The computed value function using a too small $L$ (=50, Red) and moderate $L$ (=100, Blue). Under-resolution of the nonlocal term results in a spuriously oscillating numerical solution.

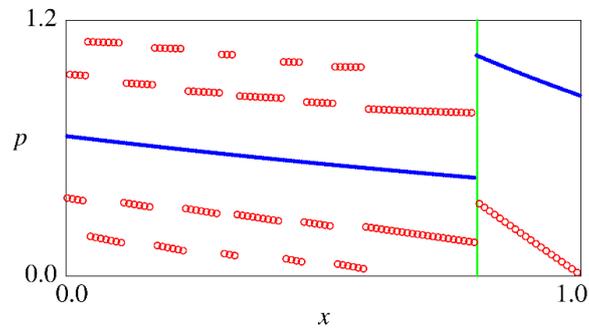

**Figure 7.** The computed PDF using a small $L$ (=100, Red) and moderate $L$ (=200, Blue). The threshold (Green) is also plotted in the figure. The numerical solution is spuriously oscillating with the small $L$. Under-resolution of the nonlocal term results in a spuriously oscillating numerical solution.

## 4.2 Application

### 4.2.1 Computational conditions

We examine impacts of designed river flows that are assumed to be specified by the environmental manager, on the coupled sediment and algae population dynamics. We consider the logistic growth $f(y) = Gy(1-y)$ with $G = 0.4$ (1/day) having the same order with the experimental and modeling results for benthic algae and periphyton [75-77]. The parameters and coefficients of the objective function and the observation process depend on the decision-makers, meaning that they are flexible. Therefore, their choice made here is only one example, and can be adaptively modeled depending on the problems and decision-makers. We assume that the decision-maker has a weekly time-scale for the system management: $\Lambda = \delta = 0.15$ (1/day). The disutility caused by the population existence is assumed to has the simplest form $S(y) = S_0 y$. We choose $S_0 = 1$. The other parameters are set as $c = 0.30$ and $d = 0.15$, which are chosen to numerically show that the threshold-type policy becomes optimal in the 2-D case as well.

The coefficient $\xi$ for the algae detachment is determined as $\xi = 4.2\bar{X}$, where $\bar{X} = 100/25$ (m$^3$/m) is the total amount of storable sediment per unit river width assuming the river width 25 (m). This parameterization is based on a physical reasoning and experimental results [59] as explained in **Appendix F**. The river flows in our model are considered in the jump process driving the system dynamics, namely the Lévy measure $v$. Assuming that the large flow events are rarer, we set $v(dz) = \lambda \chi_{\{0<z<1\}} \frac{\theta}{1-e^{-\theta}} e^{-\theta z}$ with $\theta > 0$. The parameter $\theta$ modulates the functional form of $v$: larger $\theta$ corresponds to sharper decrease of the jump size in $z$, and vice versa. The jumps larger than $z \geq 1$ is assumed to be outside the designing assuming that they are rare. In the case of the channel slope of 0.001 and the roughness coefficient of 0.03 (s/m$^{1/3}$), which are within the range of commonly found river environmental conditions, such an event is caused by the discharge around 30 (m$^3$/s) (**Appendix F**). We assume daily river discharge events and set $\lambda = 1$ (1/day). We set $\underline{z} = 0$ and $\bar{z} = 0.25$. The cut-off of $v$ at $z = \bar{z}$ is owing to its fast exponential decay. Using a larger $\bar{z}$ does not critically affect the results presented below.

The computational resolution is $n = 200$ with $L = 2n$. The pseudo-time increment for comparing the HJB equation is $\rho = 10n^{-1.5}$. The time increment for computing the FP equation is $\Delta t = n^{-1}$. The value iterations are terminated if $\left|\Phi^{(K)} - \Phi^{(K-1)}\right| < 10^{-12}$ for the HJB equation and $\left|\phi^{(k)} - \phi^{(k-1)}\right| < 10^{-10}$ ($\phi = p, q,$ and $r$) for the FP equation at all vertices.

### 4.2.2 Computational results

**Figures 8-9** show the computed value function and the optimal replenishment policy for $\theta = 50$ with the semi-Lagrangian scheme and the computed PDF using the FP equation with the FV scheme. The estimated stationary PDF using a Monte-Carlo method (10$^7$ sample paths, time increment 0.01 (day), forward Euler discretization) is provided in **Figure 10** for validating the result by the FP equation. **Figure 8** demonstrates that the optimal policy is of a threshold type (11). Comparing **Figures 9-10** shows that the computed PDFs

are close to each other between the FP equation with the FV scheme and Monte-Carlo simulation, suggesting their consistency as in the 1-D case. The sharp transition of $p$ reproduced by the Monte-Carlo method is captured in that by the FP equation. The numerical solutions do not have spurious oscillations.

The computed PDF seems to be singular along $x = 0,1$, agreeing well with the physical intuition behind the FP equation. The appearance of the singularity along $x = 0,1$ suggests the importance of carefully handling the FP equation along the boundaries. Another finding is that the PDF along the boundaries is bimodal concentrating both on large and small populations, implying that there would be clear state transitions between the large and small population states under the threshold-type optimal policy. The bimodal nature implies that the algae would thickly grow during the period of the empty state ($X_t = 0$). The appearance of the same nature during the period of fully-stored sediment state ($X_t = 1$) is because of the impulsive sediment replenishment directly from the empty state. It would be interesting to analyze in future whether this bimodal nature is observed in real cases or not.

A different value of the shape parameter representing a different flow condition is also examined as a comparison. **Figures 11-13** show the computed value function and the optimal replenishment policy for $\theta = 60$ and the computed PDFs using the corresponding FP equation with the FV scheme and the Monte-Carlo method. This is the case where the jump size is probably smaller than that with $\theta = 50$. The comparison between the cases $\theta = 50$ and 60 show that the optimal sediment replenishment policy is the threshold-type (11) in both cases, where the threshold level $\bar{x} = \bar{x}(y)$ is larger for the smaller $\theta$ with a larger probability of larger jump sizes. This means that the larger jump size more effectively suppresses the algae population as demonstrated in **Figures 9** and **12**, while it requires more frequent sediment replenishment. Furthermore, the bimodal nature of the PDF along the boundaries becomes weaker in this case, especially the boundary states for the smaller population becomes weaker, and vice versa. Again, the results between the FP equation and Monte-Carlo method are in good agreement.

Another source term $S(y) = 4\max\{y - 0.5, 0\}$ that is nonsmooth is examined with $\theta = 50$. **Figure 14** shows the computed value function and the optimal policy, suggesting that the threshold-type control policy is still optimal in this nonsmooth source case and that it is robust. The concave nature for small $y$ is weakened due to the vanishing source $S$ for $y \leq 0.5$.

Finally, the computational results of the HJB and FP equations suggest that the simplified 1-D models serve as simpler tools for analyzing river environmental management. The decision-maker can decide the sediment replenishment policy by comparing the state observations and the threshold level based on the designed flows and the objective function.

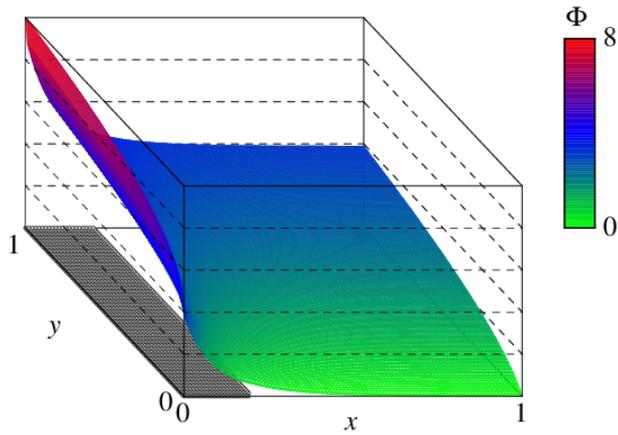

**Figure 8.** Computed value function (Surface plot) and the optimal policy (Circle plot) with $\theta = 50$.

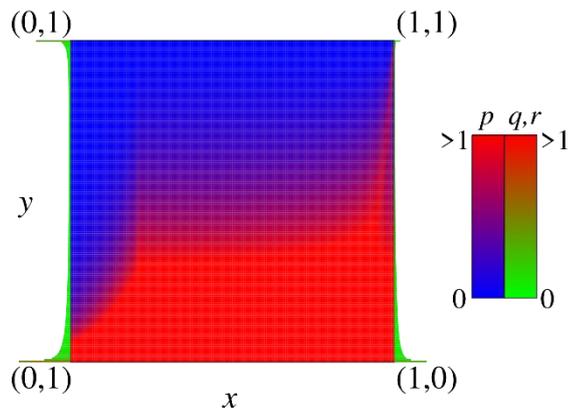

**Figure 9.** Computed PDF under the optimal policy with $\theta = 50$ by the FP equation. The maximum values of $p, q, r$ are 0.53, 0.33, 17.92, respectively.

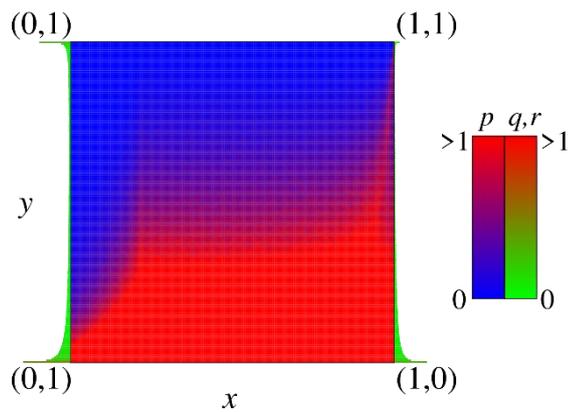

**Figure 10.** Computed PDF under the optimal policy with $\theta = 50$ by the Monte-Carlo simulation. The maximum values of $p, q, r$ are 0.49, 0.34, 18.50, respectively.

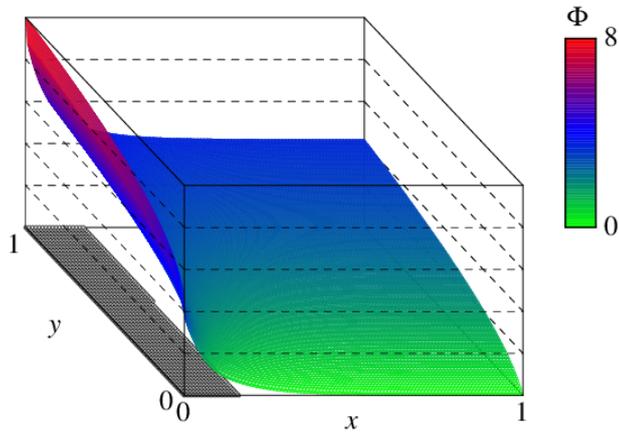

**Figure 11.** Computed value function (Surface plot) and the optimal policy (Circle plot) with $\theta = 60$.

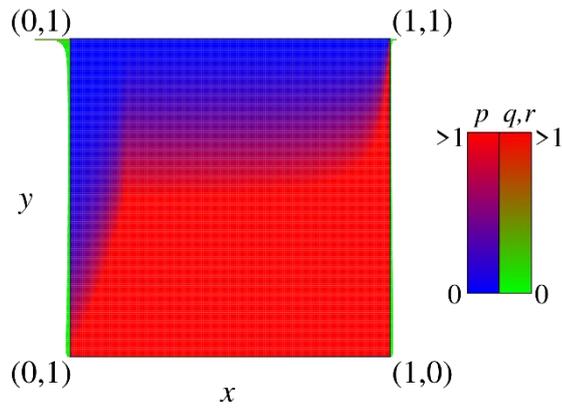

**Figure 12.** Computed PDF under the optimal policy with $\theta = 60$ by the FP equation. The maximum values of $p, q, r$ are 0.43, 0.08, 4.01, respectively.

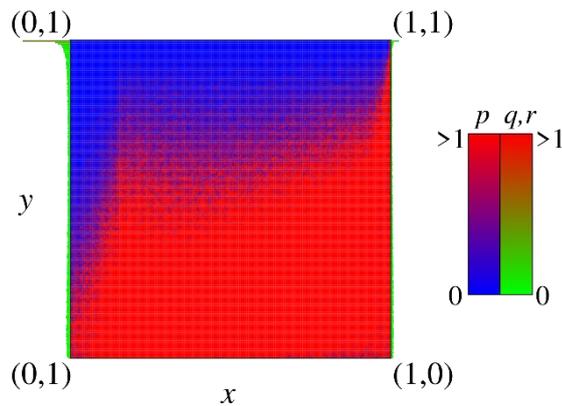

**Figure 13.** Computed PDF under the optimal policy with $\theta = 60$ by the Monte-Carlo method. The maximum values of $p, q, r$ are 0.42, 0.08, 2.42, respectively.

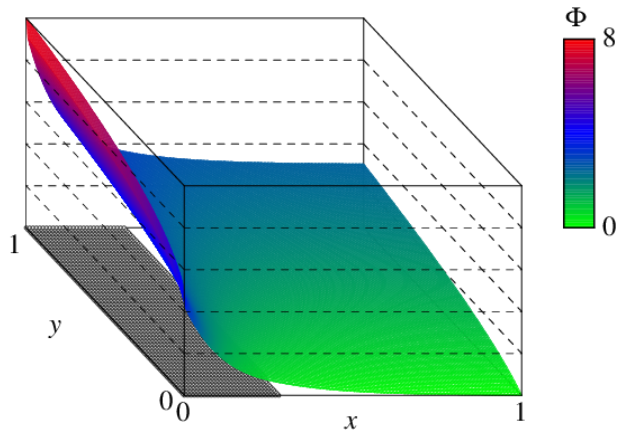

**Figure 14.** Computed value function (Surface plot) and the optimal policy (Circle plot) with a nonsmooth source term.

## 5. Conclusions

We formulated a new stochastic impulse control problem under discrete and random observation for river environmental restoration through sediment replenishment. We discussed both 2-D and 1-D cases from different viewpoints with each other, where the latter is a simplification of the former but still serving as a crucial mathematical model for better understanding mathematical and numerical discretization of our problem. The HJB and FP equations were presented and verified in the simplified case. They have boundary singularities due to the nonsmooth coefficients. Numerical schemes for handling more complicated cases were presented, demonstrating that they can correctly reproduce the singular nature of the two equations. We finally applied the presented model to demonstrative numerical computation of a more realistic case.

In this paper, only the reduced 1-D problem was mathematically analyzed, but the original coupled problem was only numerically. The HJB equation of the 2-D problem will be analyzed from a viewpoint of discontinuous viscosity solutions whose treatment would require a care especially along the discontinuity. The 2-D FP equation can be analyzed from a distributional viewpoint. Stability, monotonicity, and consistency of numerical schemes can be central properties for analyzing numerical schemes [31, 33].

We focused on a river environmental restoration problem utilizing earth and soils, but the presented mathematical framework can also be applied to other problems like the coastal erosion management by filling sand [64, 78-79] if the system dynamics are modified accordingly. In real problems, there may be some non-negligible delay between the decision-making and execution [45]. Assuming specific delay distributions, like an exponential distribution [80-81], would be an option to formulate a more advanced but still tractable model.

From a numerical viewpoint, several options to improve the presented numerical schemes would be available, which are more complicated than the presented ones. We have assumed that there is no atom along the boundaries, which would be satisfied only for not large population decay as in the simple model

[82]. In fact, the computed PDFs are in good agreement for the cases examined here. This assumption will be relaxed in the future to develop a more versatile numerical scheme that can handle singular behavior on the whole boundary.

A more sophisticated numerical scheme can be more complicated, but may be able to achieve higher-order convergence. The FP equation was found to be discontinuous even in the simplified case, implying that schemes naturally allowing for solution discontinuities, such as the discontinuous Galerkin methods [83-84], can be a more accurate alternative to our cell-centered FV scheme. For both HJB and FP equations, utilizing some moving mesh methods, like the Moving Mesh Partial Differential Equation methods [85-86], can also be an alternative to accurately resolve discontinuous solutions. For the FP equation, an energetic variational approach can be developed as a stable and convergent alternative because it is a conservation law preserving the mass. However, establishment of such a scheme is a highly non-trivial task since it requires the specification of an energy-like functional to be minimized [87-88], which has not been found for our model. Finding such a useful functional remains as a future topic. Employing a jump process with infinite activities would be both mathematical and computational challenges [89-90]. Kinetic derivation method [91] would effectively work in such a situation. We then must use a conservative implicit scheme to manage the nonlocal terms efficiently. This topic is currently under our investigations.

From an engineering viewpoint, a more comprehensive mathematical model for river environmental restoration can be established by considering dam operations [22, 92]. By using such a mathematical model, we will be able to discuss impacts of future climate changes and adaptation measures on water environmental management [93].

**Appendices**

**Appendix A: Dynamic programming**

We present a formal derivation procedure of the HJB equation (8) based on those for the existing discrete and random observations [65, 94-95]. Assume that $\tau$ is a stochastic variable following an exponential distribution with the mean $\Lambda^{-1}$. The dynamic programming principle for our problem is

$$\Phi(x,y) = \mathbb{E}^{x,y}\left[\int_0^\tau e^{-\delta s}\chi_{\{X_s=0\}}\mathrm{d}s + \int_0^\tau e^{-\delta s}S(Y_s)\mathrm{d}s + e^{-\delta \tau}\min_{\eta=\{0,1-X_\tau\}}\left\{\Phi(X_\tau+\eta,Y_\tau)+c\eta+d\chi_{\{\eta>0\}}\right\}\right]. \quad (40)$$

By the definition of $\tau$, (40) is rewritten as

$$\Phi(x,y) = \mathbb{E}^{x,y}\left[\int_0^\infty \Lambda e^{-\Lambda t}\left\{\begin{array}{l}\int_0^t e^{-\delta s}\chi_{\{X_s=0\}}\mathrm{d}s + \int_0^t e^{-\delta s}S(Y_s)\mathrm{d}s \\ +e^{-\delta t}\min_{\eta=\{0,1-X_t\}}\left\{\Phi(X_t+\eta,Y_t)+c\eta+d\chi_{\{\eta>0\}}\right\}\end{array}\right\}\mathrm{d}t\right]. \quad (41)$$

We can calculate the first integral as

$$\mathbb{E}^{x,y}\left[\int_0^\infty \Lambda e^{-\Lambda t}\int_0^t e^{-\delta s}\chi_{\{X_s=0\}}\mathrm{d}s\mathrm{d}t\right] = \mathbb{E}^{x,y}\left[\int_0^\infty \int_0^t \Lambda e^{-\Lambda t}e^{-\delta s}\chi_{\{X_s=0\}}\mathrm{d}s\mathrm{d}t\right]$$
$$= \mathbb{E}^{x,y}\left[\int_0^\infty \int_s^\infty \Lambda e^{-\Lambda t}e^{-\delta s}\chi_{\{X_s=0\}}\mathrm{d}t\mathrm{d}s\right]$$
$$= \mathbb{E}^{x,y}\left[\int_0^\infty e^{-\delta s}\chi_{\{X_s=0\}}\left(\int_s^\infty \Lambda e^{-\Lambda t}\mathrm{d}t\right)\mathrm{d}s\right], \quad (42)$$
$$= \mathbb{E}^{x,y}\left[\int_0^\infty e^{-(\delta+\Lambda)s}\chi_{\{X_s=0\}}\mathrm{d}s\right]$$

where changing the order of integrations is justified because the integrand is uniformly bounded. Similarly,

$$\mathbb{E}^{x,y}\left[\int_0^\infty \Lambda e^{-\Lambda t}\int_0^t e^{-\delta s}S(Y_s)\mathrm{d}s\mathrm{d}t\right] = \mathbb{E}^{x,y}\left[\int_0^\infty e^{-(\delta+\Lambda)s}S(Y_s)\mathrm{d}s\right]. \quad (43)$$

Consequently, we get

$$\Phi(x,y) = \mathbb{E}^{x,y}\left[\int_0^\infty e^{-(\delta+\Lambda)t}\left\{\chi_{\{X_t=0\}} + S(Y_t) + \Lambda \min_{\eta=\{0,1-X_t\}}\left\{\Phi(X_t+\eta,Y_t)+c\eta+d\chi_{\{\eta>0\}}\right\}\right\}\mathrm{d}t\right], \quad (44)$$

The right-hand side of (44) is seen as some objective function with the discount rate $\delta+\Lambda$. Applying classical Feynman-Kac formula to (44) yields (8).

**Appendix B: Proof of Proposition 1**

We show that our FP equation is conservative. The mass to be conserved is

$$M = \int_{[0,1]} q(t,y)\mathrm{d}y + \int_{[0,1]} r(t,y)\mathrm{d}y + \int_{(0,1)\times[0,1]} p(t,x,y)\mathrm{d}x\mathrm{d}y. \quad (45)$$

Assume $M(0)=1$. For $t>0$, we get

$$\frac{\mathrm{d}}{\mathrm{d}t}M = \frac{\mathrm{d}}{\mathrm{d}t}\int_{[0,1]} q(t,y)\mathrm{d}y + \frac{\mathrm{d}}{\mathrm{d}t}\int_{[0,1]} r(t,y)\mathrm{d}y + \frac{\mathrm{d}}{\mathrm{d}t}\int_{(0,1)\times[0,1]} p(t,x,y)\mathrm{d}x\mathrm{d}y$$
$$= \int_{[0,1]}\frac{\partial}{\partial t}q(t,y)\mathrm{d}y + \int_{[0,1]}\frac{\partial}{\partial t}r(t,y)\mathrm{d}y + \int_{(0,1)\times[0,1]}\frac{\partial}{\partial t}p(t,x,y)\mathrm{d}x\mathrm{d}y$$
$$= \int_{[0,1]}\left\{\begin{array}{l}-\frac{\partial}{\partial y}(f(y)q)-\Lambda q+\int_{(0,1)}\int_{(0,\infty)}\chi_{\{x<z\}}p(t,x,y)v(\mathrm{d}z)\mathrm{d}x\\+\int_{(0,\infty)}\frac{\chi_{\{y\leq g(1,z)\}}\chi_{\{z\geq 1\}}}{g(1,z)}r\left(t,\frac{y}{g(1,z)}\right)v(\mathrm{d}z)\end{array}\right\}\mathrm{d}y$$
$$+\int_{(0,1)\times[0,1]}\left\{\begin{array}{l}-\frac{\partial}{\partial y}(f(y)p)-\lambda p-\Lambda\chi_{(0,\bar{x})}p\\+\int_{(0,\infty)}\frac{\chi_{\{y\leq g(x,z)\}}\chi_{\{x<1-z\}}}{g(x,z)}p\left(t,x+z,\frac{y}{g(x,z)}\right)v(\mathrm{d}z)\\+\int_{(0,\infty)}\delta_{\{x=1-z\}}\frac{\chi_{\{y\leq g(x,z)\}}}{g(x,z)}r\left(t,\frac{y}{g(x,z)}\right)v(\mathrm{d}z)\end{array}\right\}\mathrm{d}x\mathrm{d}y$$
$$+\int_{[0,1]}\left\{-\frac{\partial}{\partial y}(f(y)r)-\lambda r+m\Lambda\right\}\mathrm{d}y$$
$$(46).$$

We calculate the equation (46) for each color.

Firstly, each integral of the red term vanishes because $f(0)=f(1)=0$. Secondly, by (17),

summing up the integrals of the green terms yields

$$\int_{[0,1]} \{-\Lambda q\} dy + \int_{(0,1)\times[0,1]} \{-\Lambda \chi_{(0,\bar{x}]} p\} dxdy + \int_{[0,1]} m\Lambda dy = \Lambda \int_{[0,1]} \left\{ -q - \int_{(0,\bar{x}]} p dx + m \right\} dy = 0. \quad (47).$$

On the blue terms, we have

$$\int_{[0,1]} \int_{(0,1]} \int_{(0,\infty)} \chi_{\{x<z\}} p(t,x,y) v(dz) dxdy = \int_{(0,\infty)} \int_{[0,1]} \int_{(0,1)} \chi_{\{x<z\}} p(t,x,y) dxdy v(dz)$$
$$= \int_{(0,\infty)} \int_{[0,1]} \int_{(0,\min\{1,z\})} p(t,x,y) dxdy v(dz) \quad (48).$$
$$= \int_{(0,\infty)} \int_{(0,\min\{1,z\})\times[0,1]} p(t,x,y) dxdy v(dz)$$

and

$$\int_{(0,1)\times[0,1]} \int_{(0,\infty)} \frac{\chi_{\{y\leq g(x,z)\}} \chi_{\{x<1-z\}}}{g(x,z)} p\left(t, x+z, \frac{y}{g(x,z)}\right) v(dz) dxdy$$
$$= \int_{(0,\infty)} \int_{(0,1)\times[0,1]} \frac{\chi_{\{y\leq g(x,z)\}} \chi_{\{x<1-z\}}}{g(x,z)} p\left(t, x+z, \frac{y}{g(x,z)}\right) dxdy v(dz)$$
$$= \int_{(0,\infty)} \int_{(0,1)\times\left[0,\frac{1}{g(x,z)}\right]} \chi_{\{u\leq 1\}} \chi_{\{x<1-z\}} p(t, x+z, u) dxdu v(dz) \quad (49).$$
$$= \int_{(0,\infty)} \int_{(0,1)\times[0,1]} \chi_{\{x<1-z\}} p(t, x+z, u) dxdu v(dz)$$
$$= \int_{(0,\infty)} \int_{(0,1)\times[0,1]} \chi_{\{x<1-z\}} p(t, x+z, y) dxdy v(dz)$$

where we utilized the transformation of variables $y = ug(x,z)$ in the third line, $0 < g \leq 1$ in the fourth line, and the trivial transformation $u = y$ in the last line. We further proceed as

$$\int_{(0,\infty)} \int_{(0,1)\times[0,1]} \chi_{\{x<1-z\}} p(t, x+z, y) dxdy v(dz)$$
$$= \int_{(0,\infty)} \int_{(0,1)\times[0,1]} \chi_{\{x+z<1\}} \chi_{\{z<1\}} p(t, x+z, y) dxdy v(dz)$$
$$= \int_{(0,\infty)} \int_{(z,1)\times[0,1]} \chi_{\{z<1\}} p(t, u, y) dudy v(dz) \quad , \quad (50)$$
$$= \int_{(0,\infty)} \int_{(\min\{1,z\},1)\times[0,1]} p(t, x, y) dxdy v(dz)$$

where we utilized the transformation of variables $x + z = u$ in the third line, and the trivial transformation $u = x$ in the last line. By (48)-(50), the sum of the integrals of the black terms becomes

$$\int_{(0,\infty)} \int_{(0,\min\{1,z\})\times[0,1]} p(t,x,y) dxdy v(dz) + \int_{(0,\infty)} \int_{(\min\{1,z\},1)\times[0,1]} p(t,x,y) dxdy v(dz)$$
$$= \int_{(0,\infty)} \int_{(0,1)\times[0,1]} p(t,x,y) dxdy v(dz) \quad . \quad (51)$$

On the brown terms, by $0 < g \leq 1$, we have

$$\int_{[0,1]}\left\{\int_{(0,\infty)}\frac{\chi_{\{y\leq g(1,z)\}}\chi_{\{z\geq 1\}}}{g(1,z)}r\left(t,\frac{y}{g(1,z)}\right)v(\mathrm{d}z)\right\}\mathrm{d}y$$

$$=\int_{[0,1]}\int_{(1,\infty)}\frac{\chi_{\{y\leq g(1,z)\}}}{g(1,z)}r\left(t,\frac{y}{g(1,z)}\right)v(\mathrm{d}z)\mathrm{d}y$$

$$=\int_{(1,\infty)}\int_{[0,1]}\frac{\chi_{\{y\leq g(1,z)\}}}{g(1,z)}r\left(t,\frac{y}{g(1,z)}\right)\mathrm{d}yv(\mathrm{d}z) \qquad , \qquad (52)$$

$$=\int_{(1,\infty)}\int_{\left[0,\frac{1}{g(1,z)}\right]}\chi_{\{u\leq 1\}}r(t,u)\mathrm{d}uv(\mathrm{d}z)$$

$$=\int_{(1,\infty)}\int_{[0,1]}r(t,u)\mathrm{d}uv(\mathrm{d}z)$$

$$=\int_{(1,\infty)}\int_{[0,1]}r(t,y)\mathrm{d}yv(\mathrm{d}z)$$

and similarly

$$\int_{(0,1)\times[0,1]}\left\{\int_{(0,\infty)}\delta_{\{x=1-z\}}\frac{\chi_{\{y\leq g(x,z)\}}}{g(x,z)}r\left(t,\frac{y}{g(x,z)}\right)v(\mathrm{d}z)\right\}\mathrm{d}x\mathrm{d}y$$

$$\int_{(0,1)\times[0,1]}\left\{\int_{(0,1)}\delta_{\{x=1-z\}}\frac{\chi_{\{y\leq g(x,z)\}}}{g(x,z)}r\left(t,\frac{y}{g(x,z)}\right)v(\mathrm{d}z)\right\}\mathrm{d}x\mathrm{d}y \qquad . \qquad (53)$$

$$=\int_{[0,1]}\int_{(0,1)}\frac{\chi_{\{y\leq g(1-z,z)\}}}{g(1-z,z)}r\left(t,\frac{y}{g(1-z,z)}\right)v(\mathrm{d}z)\mathrm{d}y$$

$$=\int_{(0,1)}\int_{[0,1]}r(t,y)\mathrm{d}yv(\mathrm{d}z)$$

By (52) and (53), the sum of the integrals of the brown terms becomes

$$\int_{(1,\infty)}\int_{[0,1]}r(t,y)\mathrm{d}yv(\mathrm{d}z)+\int_{(0,1)}\int_{[0,1]}r(t,y)\mathrm{d}yv(\mathrm{d}z)=\int_{(0,\infty)}v(\mathrm{d}z)\int_{[0,1]}r\mathrm{d}y=\lambda\int_{[0,1]}r\mathrm{d}y. \qquad (54).$$

Finally, summing up the integrals of the black terms yields

$$\int_{(0,1)\times[0,1]}\{-\lambda p\}\mathrm{d}x\mathrm{d}y+\int_{[0,1]}\{-\lambda r\}\mathrm{d}y=-\lambda\left(\int_{(0,1)\times[0,1]}p\mathrm{d}x\mathrm{d}y+\int_{[0,1]}r\mathrm{d}y\right). \qquad (55).$$

By (51), (54), and (55), the sum of the integrals of the blue, brown, and black vanishes. We can conclude that the right-hand side of (46) equals 0 and thus $\frac{\mathrm{d}}{\mathrm{d}t}M=0$, meaning that the mass is conserved.

**Appendix C: Proof of Proposition 2**

**C.1 Construction of the solution**

We seek for a non-negative and bounded solution $\overline{\Phi}\in C(0,1]$. By the assumption, we have

$$\delta\overline{\Phi}+\lambda\int_{(0,1)}\{\overline{\Phi}(x)-\overline{\Phi}(x-\min\{x,z\})\}\mathrm{d}z$$
$$+\Lambda\left(\overline{\Phi}-\min_{\eta\in\{0,1-x\}}\{\overline{\Phi}(x+\eta)+c\eta+d\chi_{\{\eta>0\}}\}\right)-\chi_{\{x=0\}}=0 \qquad \text{in }[0,1]. \qquad (56)$$

We seek for a solution such that

$$\underset{\eta\in\{0,1-x\}}{\arg\min}\left\{\bar{\Phi}(x+\eta)+c\eta+d\chi_{\{\eta>0\}}\right\}=\begin{cases}1-x & (0\le x\le \bar{x})\\ 0 & (\bar{x}<x\le 1)\end{cases}, \quad (57)$$

namely a solution such that

$$\begin{aligned}(\delta+\lambda+\Lambda)\bar{\Phi}-\lambda\int_{(0,1)}\bar{\Phi}(\max\{x-z,0\})\mathrm{d}z\\ =\Lambda\{\bar{\Phi}(1)+c(1-x)+d\}+\chi_{\{x=0\}}\end{aligned} \quad \text{for } 0\le x\le \bar{x} \quad (58)$$

and

$$(\delta+\lambda)\bar{\Phi}-\lambda\int_{(0,1)}\bar{\Phi}(\max\{x-z,0\})\mathrm{d}z=0 \quad \text{for } \bar{x}<x\le 1. \quad (59)$$

For $x>0$, we have

$$\begin{aligned}\int_{(0,1)}\bar{\Phi}(\max\{x-z,0\})\mathrm{d}z &= \int_{(0,x)}\bar{\Phi}(\max\{x-z,0\})\mathrm{d}z+\int_{(x,1)}\bar{\Phi}(\max\{x-z,0\})\mathrm{d}z\\ &= \int_{(0,x)}\bar{\Phi}(x-z)\mathrm{d}z+(1-x)\bar{\Phi}(0)\end{aligned} \quad (60)$$

and

$$\begin{aligned}\frac{\mathrm{d}}{\mathrm{d}x}\int_{(0,1)}\bar{\Phi}(\max\{x-z,0\})\mathrm{d}z &= \frac{\mathrm{d}}{\mathrm{d}x}\int_{(0,x)}\bar{\Phi}(x-z)\mathrm{d}z-\bar{\Phi}(0)\\ &= \lim_{z\to x-0}\bar{\Phi}(x-z)+\int_{(0,x)}\frac{\mathrm{d}}{\mathrm{d}x}\bar{\Phi}(x-z)\mathrm{d}z-\bar{\Phi}(0)\\ &= \bar{\Phi}(+0)-\bar{\Phi}(0)-\int_{(0,x)}\frac{\mathrm{d}}{\mathrm{d}z}\bar{\Phi}(x-z)\mathrm{d}z\\ &= \bar{\Phi}(+0)-\bar{\Phi}(0)-(\bar{\Phi}(+0)-\bar{\Phi}(x))\\ &= \bar{\Phi}(x)-\bar{\Phi}(0)\end{aligned} \quad (61)$$

By taking the limit $x\to +0$ in (58) yields

$$(\delta+\lambda+\Lambda)\bar{\Phi}(+0)-\lambda\bar{\Phi}(0)=\Lambda\{\bar{\Phi}(1)+c+d\}. \quad (62)$$

Differentiating (58) and (59) with respect to $x$ yields

$$(\delta+\lambda+\Lambda)\frac{\mathrm{d}\bar{\Phi}}{\mathrm{d}x}-\lambda(\bar{\Phi}(x)-\bar{\Phi}(0))=-\Lambda c \quad \text{for } 0<x\le \bar{x} \quad (63)$$

and

$$(\delta+\lambda)\frac{\mathrm{d}\bar{\Phi}}{\mathrm{d}x}-\lambda(\bar{\Phi}(x)-\bar{\Phi}(0))=0 \quad \text{for } \bar{x}<x\le 1, \quad (64)$$

respectively. In addition, at $x=0$, we get

$$(\delta+\Lambda)\bar{\Phi}(0)=\Lambda\{\bar{\Phi}(1)+c+d\}+1. \quad (65)$$

Solving (63) and (64) yields, with some constants $C_1, C_2$,

$$\bar{\Phi}=C_1 e^{\gamma x}+\bar{\Phi}(0)+\frac{c}{\lambda}\Lambda \quad \text{for } 0<x\le \bar{x}, \text{ with } \gamma=\lambda/(\delta+\lambda+\Lambda) \quad (66)$$

and

$$\bar{\Phi}=C_2 e^{\beta x}+\bar{\Phi}(0) \quad \text{for } \bar{x}<x\le 1, \text{ with } \beta=\lambda/(\delta+\lambda), \quad (67)$$

respectively. Taking the limit $x\to +0$ in (66) and the limit $x\to 1-0$ in (67) yields

$$C_1 = \bar{\Phi}(+0) - \bar{\Phi}(0) - \frac{c}{\lambda}\Lambda \text{ and } C_2 = -e^{-\beta}\left(\bar{\Phi}(0) - \bar{\Phi}(1)\right). \tag{68}$$

Furthermore, the continuity assumption of $\bar{\Phi}$ at $x = \bar{x}$ gives

$$\bar{\Phi}(\bar{x}) = \left(\bar{\Phi}(+0) - \bar{\Phi}(0) - \frac{c}{\lambda}\Lambda\right)e^{\gamma\bar{x}} + \bar{\Phi}(0) + \frac{c}{\lambda}\Lambda = -e^{-\beta}\left(\bar{\Phi}(0) - \bar{\Phi}(1)\right)e^{\beta\bar{x}} + \bar{\Phi}(0). \tag{69}$$

By (57) at $x = \bar{x}$, we get

$$\bar{\Phi}(\bar{x}) = \bar{\Phi}(1) + c(1 - \bar{x}) + d. \tag{70}$$

Now, we have the four unknowns $\bar{\Phi}(0), \bar{\Phi}(+0), \bar{\Phi}(1), \bar{x}$ and the four equations (62), (65), (69), (70), meaning that there are the same number of unknowns and equations. This means that if we have a solution with $\bar{x} \in (0,1)$, then $\bar{\Phi}$ satisfies the HJB equation. Finally, we rewrite $\bar{\Phi}(0), \bar{\Phi}(+0), \bar{\Phi}(1)$ as $\Phi_0, \Phi_{+0}, \Phi_1$, respectively.

**C.2 Verification**

The proof here is similar to that of Theorem 1 of Wang [65] and Proposition 2 of Yoshioka [35], but the coefficients and infinitesimal generators are different.

Set $x \in [0,1]$ and $\bar{\eta} \in C$. By Itô's formula, for any $\psi \in L^\infty([0,1]) \cap C^1((0,1])$, we get

$$e^{-\delta T}\psi(X_T) - \psi(x) = \int_0^T e^{-\delta s}\left(-\delta\psi(X_{s-}) - \int_{(0,\infty)}\{\psi(X_{s-}) - \psi(X_{s-} - \min\{X_{s-}, z\})\}\nu(\mathrm{d}z)\right)\mathrm{d}s$$
$$+ \sum_{0 \leq t \leq T} e^{-\delta t}\left(\psi(X_t + \bar{\eta}_t \mathrm{d}N_t) - \psi(X_t)\right) + M_T(\psi) \qquad , T > 0 \tag{71}$$

with a martingale term

$$M_T(\psi) = \int_0^T e^{-\delta s}\left(\int_{(0,\infty)}\{\psi(X_{s-}) - \psi(X_{s-} - \min\{X_{s-}, z\})\}\nu(\mathrm{d}z)\right)\mathrm{d}s$$
$$- \sum_{0 \leq \omega_l \leq T} e^{-\delta\omega_l}\left(\psi(X_{\omega_l -}) - \psi(X_{\omega_l})\right) \tag{72}$$

By applying (71) to $\bar{\Phi}$, we have

$$e^{-\delta T}\bar{\Phi}(X_T) - \bar{\Phi}(x) = \int_0^T e^{-\delta s}\Lambda\left(\bar{\Phi}(X_{s-}) - \min_{\eta \in \{0, 1-x\}}\{\bar{\Phi}(X_{s-} + \eta) + c\eta + d\chi_{\{\eta > 0\}}\}\right) + M_T(\bar{\Phi})$$
$$+ \int_0^T e^{-\delta s}\chi_{\{X_{s-} = 0\}}(-1)\mathrm{d}s + \sum_{0 \leq t \leq T} e^{-\delta t}\left(\bar{\Phi}(X_t + \bar{\eta}_t \mathrm{d}N_t) - \bar{\Phi}(X_t)\right) \tag{73}$$

Since $\bar{\eta}_t \neq 0$ at jumps, we have

$$\bar{\Phi}(X_t + \bar{\eta}_t \mathrm{d}N_t) - \bar{\Phi}(X_t) \geq G(X_t)\mathrm{d}N_t - \left(c + d\chi_{\{\bar{\eta}_t > 0\}}\right) \tag{74}$$

with $G:[0,1] \to \mathbb{R}$ given by

$$G(x) = \min_{\eta = \{0, 1-x\}}\{\bar{\Phi}(x + \eta) + c\eta + d\chi_{\{\eta > 0\}}\} - \bar{\Phi}(x). \tag{75}$$

Substituting (74) into (73) with the help of (75) yields

$$e^{-\delta T}\bar{\Phi}(X_T)-\bar{\Phi}(x)$$
$$=-\int_0^T e^{-\delta s}\chi_{\{X_{s-}=0\}}\,\mathrm{d}s-\int_0^T e^{-\delta s}\Lambda G(X_{s-})\,\mathrm{d}s+\sum_{0\le t\le T}e^{-\delta t}\left(\bar{\Phi}(X_t+\bar{\eta}_t\mathrm{d}N_t)-\bar{\Phi}(X_t)\right)+M_T(\bar{\Phi})$$
$$\ge -\int_0^T e^{-\delta s}\chi_{\{X_{s-}=0\}}\,\mathrm{d}s-\int_0^T e^{-\delta s}\Lambda G(X_{s-})\,\mathrm{d}s+\sum_{0\le t\le T}e^{-\delta t}\left\{G(X_t)\mathrm{d}N_t-\left(c+d\chi_{\{\bar{\eta}_t>0\}}\right)\right\}+M_T(\bar{\Phi}) \quad (76)$$
$$=-\int_0^T e^{-\delta s}\chi_{\{X_{s-}=0\}}\,\mathrm{d}s+\int_0^T e^{-\delta s}\Lambda G(X_{s-})\,\mathrm{d}\tilde{N}_s-\sum_{0\le t\le T}e^{-\delta t}\left(c+d\chi_{\{\bar{\eta}_t>0\}}\right)+M_T(\bar{\Phi})$$

with the compensated Poisson process $\tilde{N}_t=N_t-\Lambda t$ ($t\ge 0$). Taking the expectation $\mathbb{E}^x$ in both sides of (76) and rearranging it leads to

$$\bar{\Phi}(x)\le e^{-\delta T}\mathbb{E}^x\left[\bar{\Phi}(X_T)\right]+\mathbb{E}^x\left[\int_0^T e^{-\delta s}\chi_{\{X_{s-}=0\}}\,\mathrm{d}s+\sum_{0\le t\le T}e^{-\delta t}\left(c+d\chi_{\{\bar{\eta}_t>0\}}\right)\right]$$
$$-\mathbb{E}^x\left[\int_0^T e^{-\delta s}G(X_s)\,\mathrm{d}\tilde{N}_s\right] \quad (77)$$

By the Martingale property of $\tilde{N}$ and the uniform boundedness of $\bar{\Phi}$, we get

$$\mathbb{E}^x\left[\int_0^T e^{-\delta s}G(X_s)\,\mathrm{d}\tilde{N}_s\right]=0. \quad (78)$$

Taking the limit $T\to+\infty$ in (77) with the help of (78) yields

$$\bar{\Phi}(x)\le \mathbb{E}^x\left[\int_0^\infty e^{-\delta s}\chi_{\{X_s=0\}}\,\mathrm{d}s+\sum_{0\le t}e^{-\delta t}\left(c+d\chi_{\{\bar{\eta}_t>0\}}\right)\right]=\phi(x;\bar{\eta}). \quad (79)$$

Note that $X_{s-}$ in the integral can be replaced by $X_s$ because the measure $\mathrm{d}s$ is continuous. Since $\bar{\eta}\in\mathcal{C}$ is arbitrary, we get

$$\bar{\Phi}\le\Phi \quad \text{in} \quad [0,1]. \quad (80)$$

Next, we show the equality

$$\bar{\Phi}=\Phi \quad \text{in} \quad [0,1], \quad (81)$$

with which we can complete the proof. This immediately follows from the fact that the control (11) is clearly admissible. Especially, (77) becomes the equality

$$\Phi(x)=\mathbb{E}^x\left[e^{-\delta T}\Phi(X_T)+\int_0^T e^{-\delta s}\chi_{\{X_{s-}=0\}}\,\mathrm{d}s+\sum_{0\le t\le T}e^{-\delta t}\left(c+d\chi_{\{\bar{\eta}_t>0\}}\right)-\int_0^T e^{-\delta s}G(X_s)\,\mathrm{d}\tilde{N}_s\right]$$
$$=e^{-\delta T}\mathbb{E}^x\left[\Phi(X_T)\right]+\mathbb{E}^x\left[\int_0^T e^{-\delta s}\chi_{\{X_{s-}=0\}}\,\mathrm{d}s+\sum_{0\le t\le T}e^{-\delta t}\left(c+d\chi_{\{\bar{\eta}_t>0\}}\right)\right] \quad (82)$$

with this control by (78). We again take the expectation $\mathbb{E}^x$ and the limit $T\to+\infty$ and obtain (81).

**Appendix D: Proof of Proposition 3**

Firstly, assume $\bar{x}>0$. The mass conservation condition gives

$$q+r+\int_{(0,1)}p(x)\mathrm{d}x=1. \quad (83)$$

The equation (21) is rewritten as

$$(\lambda+\Lambda)p=\lambda\int_{(0,1-x)}p(x+z)\mathrm{d}z+\lambda r=\lambda\int_{(x,1)}p(u)\mathrm{d}u+\lambda r,\quad 0<x\le\bar{x} \quad (84)$$

with

$$\lambda p = \lambda \int_{(x,1)} p(u) du + \lambda r, \quad \bar{x} < x < 1. \tag{85}$$

The equation (85) is directly solved as

$$p(x) = re^{1-x}, \quad \bar{x} < x < 1. \tag{86}$$

Differentiating (84) with respect to $x$ gives

$$p(x) = C_1 e^{-\alpha x}, \quad 0 < x \leq \bar{x} \tag{87}$$

with $\alpha = \lambda/(\lambda+\Lambda)$ and a constant $C_1$. Substituting $x = \bar{x}$ into (84) with (86) gives $C_1 = \alpha r e^{1-\bar{x}+\alpha\bar{x}}$.
We then have

$$\begin{aligned}
q &= \lambda\Lambda^{-1} r - \int_{(0,\bar{x})} p(x) dx \\
&= r\left(\lambda\Lambda^{-1} - \int_{(0,\bar{x})} \alpha r e^{1-\bar{x}+\alpha\bar{x}} e^{-\alpha x} dx\right) \\
&= r\left(\lambda\Lambda^{-1} - e^{1-\bar{x}}\left(e^{\alpha\bar{x}} - 1\right)\right)
\end{aligned} \tag{88}$$

by the second equation of (22). From (83), we get

$$\begin{aligned}
1 &= q + r + \int_{(0,1)} p(x) dx \\
&= q + r + \int_{(0,\bar{x})} p(x) dx + \int_{(\bar{x},1)} p(x) dx \\
&= r(1 + \lambda\Lambda^{-1}) + \int_{(\bar{x},1)} p(x) dx \\
&= r(1 + \lambda\Lambda^{-1}) + r(e^{1-\bar{x}} - 1)
\end{aligned} \tag{89}$$

and thus $r = \left(\lambda\Lambda^{-1} + e^{1-\bar{x}}\right)^{-1} > 0$, from which we get $q$ and $p$ as well. Notice that we have not used the first equation of (22). However, it is actually a redundant and is automatically satisfied by the mass conservation (83) and the fact that $q$ is a scalar.

Finally, $q > 0$ follows from the fact that the quantity $\lambda\Lambda^{-1} - e^{1-\bar{x}}\left(e^{\alpha\bar{x}} - 1\right) = u + e^{1-\bar{x}} - e^{1-\frac{1}{1+u}\bar{x}} \left(\equiv F(\bar{x})\right)$ with $u = \lambda\Lambda^{-1} > 0$, as a function of $\bar{x} \in [0,1]$, is decreasing and satisfies $F(1) = 1 + u - e^{1-\frac{1}{1+u}} = (1+u)e^{-\frac{1}{1+u}}\left(e^{\frac{1}{1+u}} - \frac{1}{1+u}e\right) > 0$. Note the elementary result $e^\xi - e\xi \geq 0$ ($\xi \geq 0$).

## Appendix E: Proof of Proposition 4

It is sufficient to show $\dfrac{dM}{dt} = 0$ for $t > 0$. Fix $t > 0$. We have

$$\begin{aligned}
\frac{\mathrm{d}M}{\mathrm{d}t} &= \sum_{j=1}^{n}\frac{\mathrm{d}q_j}{\mathrm{d}t}h + \sum_{i,j=1}^{n}\frac{\mathrm{d}p_{i,j}}{\mathrm{d}t}h^2 + \sum_{j=1}^{n}\frac{\mathrm{d}r_j}{\mathrm{d}t}h \\
&= -\sum_{j=1}^{n}\left(\Lambda q_j + \frac{1}{h}\left(F_{q,j} - F_{q,j-1}\right) - J_j^{(L)}\right)h \\
&\quad -\sum_{i,j=1}^{n}\left(\Lambda \chi_{\{i\le\theta_j\}} p_{i,j} + \frac{1}{h}\left(F_{p,i,j} - F_{p,i,j-1}\right) + \lambda p_{i,j} - J_{i,j}\right)h^2 \\
&\quad -\sum_{j=1}^{n}\left(\frac{1}{h}\left(F_{r,j} - F_{r,j-1}\right) + \lambda r_j - \Lambda m_j\right)h
\end{aligned} \qquad (90)$$

The conservation property of the flux terms containing $F_{q,i,j}$, $F_{p,i,j}$, and $F_{r,j}$ are easy to check by their definitions and the vanishing property along the boundaries. In addition, we have

$$-\sum_{j=1}^{n} q_j h - \sum_{i,j=1}^{n}\Lambda\chi_{\{i\le\theta_j\}} p_{i,j} h^2 - \sum_{j=1}^{n}(-m_j)h = -\sum_{j=1}^{n} m_j h + \sum_{j=1}^{n} m_j h = 0 \qquad (91)$$

by the definition of $m_j$. Then, (90) reduces to

$$\frac{\mathrm{d}M}{\mathrm{d}t} = \left(\sum_{j=1}^{n} J_j^{(L)} h + \sum_{i,j=1}^{n} J_{i,j} h^2\right) - \lambda\left(\sum_{i,j=1}^{n} p_{i,j} h^2 + \sum_{j=1}^{n} r_j h\right). \qquad (92)$$

The first term in the right-hand side of (92) is calculated as

$$J_{i,j} = \sum_{l=1}^{L}\sum_{i',j'=1}^{n}\chi_{\{\alpha_{i',l}=i-1,\beta_{i',j',l}=j-1\}} v_l p_{i',j'} + \sum_{l=1}^{L}\sum_{j'=1}^{n}\chi_{\{\gamma_l=i-1,\omega_{j',l}=j-1\}} v_l r_{j'} h^{-1}, \qquad (93)$$

$$J_j^{(L)} = \sum_{l=1}^{L}\sum_{i',j'=1}^{n}\chi_{\{\alpha_{i',l}<0,\beta_{i',j',l}=j-1\}} p_{i',j'} h + \sum_{l=1}^{L}\sum_{j'=1}^{n}\chi_{\{\gamma_l<0,\omega_{j',l}=j-1\}} r_{j'}. \qquad (94)$$

$$\begin{aligned}
\sum_{j=1}^{n} J_j^{(L)} h + \sum_{i,j=1}^{n} J_{i,j} h^2 &= h\sum_{j=1}^{n}\sum_{l=1}^{L}\sum_{i',j'=1}^{n}\chi_{\{\alpha_{i',l}<0,\beta_{i',j',l}=j-1\}} v_l p_{i',j'} h + h\sum_{j=1}^{n}\sum_{l=1}^{L}\sum_{j'=1}^{n}\chi_{\{\gamma_l<0,\omega_{j',l}=j-1\}} v_l r_{j'} \\
&\quad + h^2\sum_{i,j=1}^{n}\sum_{l=1}^{L}\sum_{i',j'=1}^{n}\chi_{\{\alpha_{i',l}=i-1,\beta_{i',j',l}=j-1\}} v_l p_{i',j'} + h^2\sum_{i,j=1}^{n}\sum_{l=1}^{L}\sum_{j'=1}^{n}\chi_{\{\gamma_l=i-1,\omega_{j',l}=j-1\}} v_l r_{j'} h^{-1} \\
&= h^2\sum_{j=1}^{n}\sum_{l=1}^{L}\sum_{i',j'=1}^{n}\chi_{\{\beta_{i',j',l}=j-1\}} v_l p_{i',j'} + h\sum_{j=1}^{n}\sum_{l=1}^{L}\sum_{j'=1}^{n}\chi_{\{\omega_{j',l}=j-1\}} v_l r_{j'} \\
&= h^2\sum_{i',j'=1}^{n} p_{i',j'}\sum_{l=1}^{L} v_l \sum_{j=1}^{n}\chi_{\{\beta_{i',j',l}=j-1\}} + h\sum_{j'=1}^{n} r_{j'}\sum_{l=1}^{L} v_l \sum_{j=1}^{n}\chi_{\{\omega_{j',l}=j-1\}}
\end{aligned} \qquad (95)$$

We have used the fact that $\sum_{i=1}^{n}\chi_{\{\alpha_{i',l}=i-1\}} + \chi_{\{\alpha_{i',l}<0\}} = 1$ for each couple $(i',l)$ because $\alpha_{i',l}$ is either negative or exactly one of $\{0,1,2,...,n-1\}$. For each triplet $(i',j',l)$, there is a unique natural number $1\le j\le n$ with $\chi_{\{\beta_{i',j',l}=j-1\}}=1$, meaning $\sum_{j=1}^{n}\chi_{\{\beta_{i',j',l}=j-1\}}=1$. Similarly, $\sum_{j=1}^{n}\chi_{\{\omega_{j',l}=j-1\}}=1$. By the definition, $\sum_{l=1}^{L} v_l = \lambda$. Consequently,

$$\sum_{j=1}^{n} J_j^{(L)} h + \sum_{i,j=1}^{n} J_{i,j} h^2 = h^2 \sum_{i',j'=1}^{n} p_{i',j'} \sum_{l=1}^{L} v_l + h \sum_{j'=1}^{n} r_{j'} \sum_{l=1}^{L} v_l = \lambda \left( \sum_{i',j'=1}^{n} p_{i',j'} h^2 + \sum_{j'=1}^{n} r_{j'} h \right). \tag{96}$$

Combining (92) and (96) yields the desired result, and completes the proof.

**Appendix F: Physical reasoning of the algae detachment modeling**

We give a physical reasoning of the coefficient $g$ for the algae detachment. It has been experimentally found that the algae detachment due to river flows containing sediment (bedload) particles occurs within 1 to several hours [59], which is much shorter than the time-scale of algae growth as it is daily to weekly. The exponential functional form is motivated by the same experimental result that the algae population decreases exponentially in time. Assume enough sediment storage unless otherwise specified.

The sediment transport per unit time and width is denoted as $Q$, and the duration of a high flood during which the algae detachment occurs as $h$. The total amount of the sediment transported toward downstream during this period is $Qh$. The experimentally found exponential decrease of the algae population can be expressed as $Y_t = Y_{t-h} \exp(-\mu h)$ with some $\mu > 0$ depending on the flow and sediment transport conditions [59, 96]. This $\mu$ has been found to be reasonably approximated as a linear function of the transport rate $Q$: $\mu = \kappa Q$ with some $\kappa > 0$. Consequently, we get $Y_t = Y_{t-h} \exp(-\kappa Qh) = Y_{t-h} \exp(-\kappa(X_{t-h} - X_t))$ because of $X_{t-h} - X_t = Qh$. By normalizing $X$ with the storable amount of the sediment $\bar{X} > 0$ yields $Y_t = Y_{t-h} \exp(-\kappa \bar{X}(X_{t-h} - X_t))$, where $X$ here is non-dimensional and in $[0,1]$. Then, formally taking the limit $h \to 0$ yields $Y_t = Y_{t-} \exp(-\kappa \bar{X} z_t)$, where $z_t$ is the jump size.

The discussion above assumed that the sediment storage is sufficiently large, but we can effectively replace $z_t$ by $\min\{z_t, X_{t-}\}$. We then get $Y_t = Y_{t-} \exp(-\kappa \bar{X} z_t)$, meaning that we should choose $\xi = \kappa \bar{X}$. In this way, we have to specify the two parameters: the amount of storable sediment $\bar{X}$ and the regression parameter $\kappa$ whose value can be found experimentally. From the experimental result with the sand particle having the diameter of 0.005 (m) [59], we get $\kappa = 4.2$ (1/m²). Of course, a more complicated model in Yoshioka et al. [47] can also be used but the parameter values will be more sensitive.

By using the Meyer-Peter-Muller and Manning's formulae [35, 49] for the particle diameter 5 (mm), river width of 25 (m), and the Manning's roughness coefficient of 0.03 (s/m$^{1/3}$), the bed slope of 0.001 yields the formula relating the unit-width sediment discharge $Q$ and water flow discharge as $Q = 0.0112 \max\{0.0176 Q_w^{0.6} - 0.047, 0\}^{1.5}$ (m²/s) with $Q_w$ (m³/s) the water flow discharge. This formula shows that the sediment transport occurs only when the flow discharge is sufficiently large with which the sediment particles can be transported along the river. Then, we get $X_{t-h} - X_t = Q(Q_w)h$ whose right-hand side has a function $Q(Q_w)$ of the flow discharge $Q_w$. This means $z_t = Q(Q_w)h$ provided that

$X_{t-h}$ is sufficiently large. By the above-specified parameters, setting $h = O(10^0)$ (h), the time-scale of detachment, with $z = 1$ corresponds to $Q_w = O(10^2)$ (m³/s), which is a sufficiently high discharge such that the sediment is completely flushed out by the river flows.

**Acknowledgements**

JSPS Research Grant 18K01714 and 19H03073, and grants from MLIT Japan for ecological survey of a life history of the landlocked *Plecoglossus altivelis altivelis* and management of seaweed in Lake Shinji support this research. The explanation of our mathematical model was significantly improved by valuable comments from an anonymous reviewer.

**Declarations of interest**

None.